\documentclass[10pt,fleqn]{amsart}

\usepackage[utf8]{inputenc}
\usepackage{amsmath,amssymb,latexsym}
\usepackage[mathscr]{eucal}
\usepackage{accents}
\usepackage[dvipsnames]{xcolor}

\definecolor{labelkey}{gray}{0.80}

\usepackage{cite}
\usepackage{psfrag}

\usepackage{eqnnumwarn}

\theoremstyle{plain}
\newtheorem{theorem}{Theorem}
\newtheorem{prop}[theorem]{Proposition}
\newtheorem{lemma}{Lemma}

\numberwithin{equation}{section}

\def\mydate{\number\year-\ifnum\month<10{0}\fi\number\month-\ifnum\day<10{0}\fi\number\day}

\newcounter{mnote}
\newcommand{\dy}{\partial}

\newcommand{\sfrac}[2]{{\textstyle\frac{#1}{#2}}}

\newcommand{\tssum}{{\textstyle\sum}}
\newcommand{\tssuml}{\textstyle\sum\limits} 

\newcommand{\Zahl}{\mathbb{Z}}
\newcommand{\Real}{\mathbb{R}}

\newcommand{\ex}{\mathrm{e}}
\newcommand{\im}{\mathrm{i}}
\newcommand{\eps}{\varepsilon}
\newcommand{\vfi}{\varphi}
\newcommand{\vtt}{\vartheta}

\newcommand{\gb}{\nabla}

\newcommand{\gbi}{|\nabla|^{-1}}
\newcommand{\agb}{|\nabla|}

\newcommand{\aand}{\quad\textrm{and}\quad}

\newcommand{\Rar}{\quad\Rightarrow\quad}

\newcommand{\Dom}{D}

\newcommand{\kpg}{\kappa_g}
\newcommand{\kpb}{\bar\kappa}
\newcommand{\kpe}{\kappa_\eta}
\newcommand{\Proj}{{\sf P}}
\newcommand{\Zp}{\mathbb{Z}^2\backslash\{0\}}
\newcommand{\Zu}{\mathbb{Z}_+^2}

\newcommand{\uh}{\hat u}
\newcommand{\vh}{\hat v}
\newcommand{\vtth}{\hat\vtt}
\newcommand{\mb}{{\nu}} 
\newcommand{\tht}{\theta}
\newcommand{\tco}{\Theta} 
\newcommand{\tcoh}{\hat\tco}
\newcommand{\thth}{\hat\theta}
\newcommand{\thts}{\vartheta} 

\newcommand{\dvtt}{\delta\vtt}
\newcommand{\vttH}{\tilde\vtt}

\newcommand{\dr}{\;\mathrm{d}r}
\newcommand{\ilap}{\Delta^{-1}}

\newcommand{\Unif}{\mathcal{U}}
\newcommand{\EV}{{\sf E}\mkern 2mu}

\begin{document}

\title[BHT spectrum with large velocity]
{The Batchelor--Howells--Townsend spectrum:\\large velocity case}
\author{M. S. Jolly$^{1}$}
\address{$^1$Mathematics Department\\
Indiana University\\ Bloomington, IN 47405, United States}
\author{D. Wirosoetisno$^{2}$}
\address{$^2$Mathematical Sciences\\
Durham University\\ Durham, United Kingdom\ \ DH1 3LE}
\email[M. S. Jolly]{msjolly@indiana.edu}
\email[D. Wirosoetisno]{djoko.wirosoetisno@durham.ac.uk}

\thanks{This work was supported in part by Simons Foundation grant  MP-TSM-00002337.
We are grateful to C.~Demeter for the proof of \eqref{q:riesz}.}


\begin{abstract}
We consider the behaviour of a passive tracer $\tht$ governed by $\dy_t\tht + u\cdot\gb\tht = \Delta\tht + g$ in two space dimensions with prescribed smooth random incompressible velocity $u(x,t)$ and source $g(x)$.
In 1959, Batchelor, Howells and Townsend (J.\ Fluid Mech.\ 5:113) predicted that the tracer (power) spectrum should then scale as $|\tht_k|^2\propto|k|^{-4}|u_k|^2$ for $|k|$ large depending on the velocity $u$.
For smaller $|k|$, Obukhov and Corrsin earlier predicted a different spectral scaling.
In this paper, we prove that the BHT scaling does indeed hold probabilistically for sufficiently large $|k|$, asymptotically up to controlled remainders, using only bounds on the smaller $|k|$ component.
\end{abstract}

\maketitle


\section{Introduction}\label{s:intro}

The study of statistical behaviour of passive tracers in turbulent flows is nearly as old as the statistical theory of turbulence itself: eight years after A.N.\ Kolmogorov's seminal work, his student A.M.\ Obukhov \cite{obukhov:49} proposed that in a three-dimensional homogeneous isotropic turbulent flow with energy spectrum scaling as $E(k)\propto k^{-5/3}$, the power spectrum of a passive tracer should also scale as $T(k)\propto k^{-5/3}$, with the same prediction made almost contemporaneously by Corrsin \cite{corrsin:51}.
As there are several parameters involved, namely the transfer rates of energy and of tracer (variance), viscosity (of the fluid) and diffusivity (of the tracer), the full picture is more complicated.
It was soon clarified that the above Obukhov--Corrsin regime is to hold only in the energy inertial range (between the forcing and Kolmogorov scales) and where the tracer diffusivity is weak compared to the velocity.
In those parts of the inertial range where tracer diffusivity predominates, Batchelor, Howard and Townsend \cite[henceforth BHT]{batchelor-howells-townsend:59} argued that the tracer spectrum should scale as $T(k)\propto k^{-4}E(k)$, so $k^{-17/3}$ in 3d.
In the dissipation range (smaller than the Kolmogorov scale), one has the Batchelor $k^{-1}$ spectrum \cite{batchelor:59} when the diffusivity is still negligible.
Both energy and tracer spectra decay exponentially in the smallest scales.

A similar picture holds in two space dimensions, where Kraichnan \cite{kraichnan:67}, Leith \cite{leith:68} and Batchelor \cite{batchelor:69} proposed an enstrophy inertial range where $E(k)\propto k^{-3}$.
Due to intermittency, logarithmic corrections to this have been proposed by various authors.
At the smallest scales, beyond the ``Kraichnan scale'', the energy spectrum is again expected to decay exponentially.
As for the tracer spectrum in the enstrophy inertial range, analogous argument gives $T(k)\propto k^{-1}$ where diffusivity is weak (OC spectrum), and $T(k)\propto k^{-7}$ where diffusivity predominates (BHT spectrum).
For modern synopses of these physical arguments, see, e.g., \cite{davidson:t,lesieur:tf,vallis:aofd} for turbulence generally, and \cite{gotoh-yeung:13,holzer-siggia:94,shraiman-siggia:00,warhaft:00,sreenivasan:19} for tracers specifically;
we also refer to \cite{rubinstein-clark:12} for closure approaches to the tracer problem, noting in particular that the ``classical'' closures do not produce the BHT spectrum.

While invaluable for their important dynamical insight and deep physical intuition, these arguments made many assumptions still awaiting mathematical justifications.
Leaving aside the much more difficult case of the Navier--Stokes equations (in both two and three dimensions), there remain many open aspects of the passive tracer spectrum problem.
An early breakthrough was obtained by Kraichnan \cite{kraichnan:68}, using a rough random velocity whose correlation timescale is less than the diffusive one, allowing the problem to be solved exactly in all the Obukhov--Corrsin, BHT and Batchelor regimes.
There have been many related works following this paper, notably \cite{kraichnan:94,gawedzki-kupiainen:95}; we refer to \cite{eyink-jafari:23} for references.
As noted by Kraichnan himself \cite[p.~947, para.~2]{kraichnan:68}, his setup uses ``assumptions markedly different from Batchelor's''; we believe that one should therefore consider BHT59 and K68 as treating different problems.

More recently, \cite{bedrossian-blumenthal-punshon:22} rigorously established the Batchelor spectrum with velocity being a solution of stochastic Navier--Stokes equations (with hyperviscosity in 3d), so much smoother than K68, and tracer source that is smooth in space but white in time.

It has been pointed out \cite{sreenivasan:19,gotoh-yeung:13} that, in contrast with the Kolmogorov $k^{-5/3}$ or Kraichnan $k^{-3}$ spectra, the BHT spectrum is difficult to observe in nature (the only known material allowing very small Schmidt numbers and large Reynolds numbers are liquid metals) and in numerical simulations with realistic velocity fields (since one is to resolve both the OC and BHT regimes).

Following our previous works \cite{jolly-dw:bht,jolly-dw:b3d}, in this paper we consider BHT59's original problem, with a synthetic random velocity that is smooth in time and space (but having independent Fourier modes).
Unlike in \cite{jolly-dw:bht,jolly-dw:b3d}, here we make no assumptions on the size of the velocity, so our tracer spectrum will in general consist of both a (presumed) Obukhov--Corrsin part and a BHT part.
Being unable to handle the OC part, we ``model'' the lower modes by a pseudospectral system $\tco$ and prove that the higher modes have the BHT scaling \eqref{q:evtt-ukpe}, both up to controlled remainders.

One limitation of the current approach is that it cannot handle energy spectrum shallower than $k^{-3}$ in two dimensions; in fact, this K67 \cite{kraichnan:67} energy spectrum is the marginal case for our approach, necessitating an ultraviolet cutoff.
It appears to us that new mathematical ingredients would be needed to extend this type of result to three dimensions with the Kolmogorov $k^{-5/3}$ energy spectrum.
Another limitation is that we are only able to establish \eqref{q:evtt-ukpe} for wavenumbers $|k|\gtrsim c\,\|u\|_\infty^2$, whereas the classical prediction gives the transition wavenumber to be $\propto U^{1/2}$; this is discussed further following Thm.~\ref{t:uht-ukpe}.

\bigskip\hbox to\hsize{\qquad\hrulefill\qquad}\medskip

We consider the evolution of a passive scalar $\tht(x,t)$ under a prescribed incompressible velocity field $u(x,t)$ and source $g(x)$,
\begin{equation}\label{q:dtht}
   \dy_t\tht + u\cdot\gb\tht = \Delta\tht + g.
\end{equation}
For simplicity, we take $x\in\Dom:=[0,2\pi]^2$ and assume periodic boundary conditions in all directions (sometimes we will write $d$ instead of $2$ to emphasize the role of the dimension).

As noted in the original BHT paper \cite{batchelor-howells-townsend:59} and justified explicitly in \cite{jolly-dw:bht}, the case with time-dependent velocity $u(x,t)$ is not fundamentally different from the static case $u(x)$ provided that the time derivative $\dy_t u$ is sufficiently regular.
This is in sharp contrast with the case of noisy-in-time $u$ treated in \cite{kraichnan:68} and its many subsequent works, where the weak temporal correlations of $u$ is crucial in their argument; as we argued above, this should be considered as a different problem.
In this paper, therefore, we shall consider static velocity $u(x)$, so here $\tht$ is the solution of
\begin{equation}\label{q:tht}
   -\Delta\tht + u\cdot\gb\tht = g.
\end{equation}
With no loss of generality, we assume that
\begin{equation}\label{q:mean0}
   \int_\Dom u(x) \;\mathrm{d}x = 0
   \quad\text{and}\quad
   \int_\Dom \tht(x) \;\mathrm{d}x = 0.
\end{equation}
For consistency, we must impose the same condition on $g$. 
Due to the nature of the problem, all functions have zero integral over $\Dom$; therefore all our $L^p$ and Sobolev $H^s$ spaces are understood to be homogeneous.

Assuming that $\tht$ (and $u$ later) is at least $L^2$, we expand $\tht(x)$ in Fourier series
\begin{equation}\label{q:four}
  \tht(x) = \tssum_{k\in\Zahl^2}\,\thth_k\ex^{\im k\cdot x},
\end{equation}
but \eqref{q:mean0} implies that $\thth_0=0$;
throughout this paper, sums over wavenumbers always exclude $k=0$.
For $1\le\kappa<\kappa'$, we denote spectral projection by
\begin{equation}\label{q:Pdef}
   (\Proj_{\kappa,\kappa'}\tht)(x) := \tssum_{\kappa\le|k|<\kappa'} \thth_k\ex^{\im k\cdot x},
\end{equation}
and write for conciseness (our $\kappa$-intervals are always closed below and open above),
\begin{equation}
   v^{<\kappa} := \Proj_{1,\kappa}v,
   \quad
   v^{>\kappa} := \Proj_{\kappa,\infty}v
   \aand
   v_{\kappa,\kappa'} := \Proj_{\kappa,\kappa'}v.
\end{equation}
To keep the proofs readable, we assume that $g$ is spectrally bounded,
\begin{equation}\label{q:g}
   g(x) = \tssum_{1\le|k|<\kpg}\, \hat g_k \ex^{\im k\cdot x}
\end{equation}
for some $\kpg<\infty$; we expect our results to extend to the more general case of sufficiently smooth $g$.

\noindent {\bf Notation:}
We denote by $\|v\|_p^{}$ the (Lebesgue measure) $L^p$ norm in $[0,2\pi)^2$, and similarly by $\|\hat v\|_p^{}$ the (counting measure) $\ell^p$ norm in $\Zp$.
From the definition \eqref{q:four} and Parseval/Plancherel, we have (these will often be used with no mention)
\begin{equation}
   \|\tht\|_2^2 = (2\pi)^2\|\thth\|_2^2
   \aand
   \|\tht\|_\infty^{} \le \|\thth\|_1^{}.
\end{equation}
When the meaning is clear, we often write
$(v,w):=(v,w)_{L^2}^{}$ for brevity.
We denote by $*$ both the continuous and discrete convolutions, thus
\begin{equation}
   (f*g)(x) = \int_\Dom f(x-y)g(y) \;\mathrm{d}y
   \aand
   (\hat f * \hat g)_k = \tssum_j\, \hat f_{k-j} \hat g_j.
\end{equation}
They both satisfy Young's inequality: for $1\le p,q,r \le \infty$ with $1+1/r=1/p+1/q$,
\begin{equation}\label{q:young}
   \|f*g\|_r^{} \le \|f\|_p^{}\|g\|_q^{}
   \aand
   \|\hat f*\hat g\|_r^{} \le \|\hat f\|_p^{}\|\hat g\|_q^{}.
\end{equation}

For any $s\in\Real$, and for real-valued $f$ implicitly having zero integral, we define $\agb^s f:=(-\Delta)^{s/2} f$ formally by Fourier series,
\begin{equation}
   (\agb^s f)(x) := \tssum_{k\ne0}\,|k|^s \hat f_k \ex^{\im k\cdot x}
\end{equation}
and similarly component-wise for $\agb^s u$.
From the definition, $\|\agb f\|_2^{}=\|\gb f\|_2^{}$ so we will usually use the latter expression for short; however $\|\agb f\|_\infty^{}\not\sim_c \|\gb f\|_\infty^{}$.
We note some useful inequalities (see the Appendix for proofs):
\begin{align}
   &\|\gbi\gb\cdot u\|_2^{} \le c\,\|u\|_2^{} \label{q:gbiu}\\
   &\|\gbi\gb\cdot u\|_\infty^{} \le c\,\|\uh\|_1^{}
	\le c\,\|u\|_2^{1/2}\|\Delta u\|_2^{1/2} \label{q:agmon}\\
   &\|\gbi\gb\cdot u\|_\infty^{} \le c\,\|\uh\|_1^{}
	\le c\,\|\gb u\|_2^{}\bigl(1 + \log(\|\Delta u\|_2^{}/\|\gb u\|_2^{})^{1/2}\bigr). \label{q:bgu}
\end{align}
The last two inequalities show that $\|\gbi\gb\cdot u\|_\infty^{}$ is bounded by, respectively, the Agmon and Brezis--Gallou{\"e}t inequalities for $\|u\|_\infty^{}$.
Bounding $\|\gbi\gb h\|_\infty^{}\le\|\hat h\|_1^{}$ is rather wasteful, and we will need a sharper bound for later use.

For $m\in\{0,1,2,\cdots\}$ and $r>0$, we define $L_m(\cdot,\cdot)$ by
\begin{equation}\begin{aligned}
   \int_\kappa^\infty &k^{-r-1}\log^mk \;\mathrm{d}k
	= \frac{\kappa^{-r}\log^m\kappa}{r} + \frac{m}{r} \int_\kappa^\infty k^{-r-1}\log^{m-1}k \;\mathrm{d}k\\
	&= \kappa^{-r}r^{-1}\bigl( \log^m\kappa + m/r\log^{m-1}\kappa + \cdots + m!/r^{m}\bigr)
	=: L_m(\kappa,r)/(r\kappa^r).
\end{aligned}\end{equation}
As we may have to consider the case of small $s$, in this section $c(s)$ denotes a generic constant depending on $s$ that is bounded as $s\to0^+$ while the generic $c$ does not depend on $s$ (or any other parameter).

\begin{lemma}\label{t:riesz}
For $s>0$, $p\in[1,\infty]$ and $w\in L^p(\Dom)$ with zero integral, we have
\begin{align}
   &\|\agb^{-s}w^{>\kappa}\|_p^{} \le c(s)\kappa^{-s}L_2(2\kappa,s) \|w^{>\kappa}\|_p^{} & &\text{for }\kappa\ge5, \label{q:logb}\\
   &\|\agb^{-s}w\|_p^{} \le c(s)s^{-2}\|w\|_p^{}, \label{q:bern1}\\
   &\|\agb^{-s-1}\gb w^{>\kappa}\|_p^{} \le c(s)s^{-1}\kappa^{-s}L_2(2\kappa,s) \|w^{>\kappa}\|_p^{} & &\text{for }\kappa\ge5, \label{q:riesz-s}\\
   &\|\agb^{-s-1}\gb w\|_p^{} \le c(s)s^{-3}\|w\|_p^{}. \label{q:riesz}
\end{align}
\end{lemma}

\noindent We note that this result is sub-optimal for $p<\infty$, since (the case $s=0$) the Riesz transform $\agb^{-1}\gb$ is known to be bounded in $L^p$ whereas our bound \eqref{q:riesz} is infinite, but it suffices for our purposes here.
The proof can be found in Appendix~\ref{s:riesz}.

We collect some useful identities and inequalities involving $\Proj$.
The ideas are well-known, but for completeness we prove \eqref{q:poi}--\eqref{q:kinfty} in the Appendix.
\begin{lemma}\label{t:P}
Given $\kappa\ge1$, for any zero-mean $v$ and $w$ we have, whenever the expressions are defined,
\begin{align}
  &(v,w^{>\kappa})_{L^2}^{} = (v^{>\kappa},w^{>\kappa})_{L^2}^{} = (v^{>\kappa},w)_{L^2}^{}\\
  &(v,w^{<\kappa})_{L^2}^{} = (v^{<\kappa},w^{<\kappa})_{L^2}^{} = (v^{<\kappa},w)_{L^2}^{}\\
  &\|v^{>\kappa}\|_{2}^{} \le \kappa^{-1}\,\|\gb v^{>\kappa}\|_{2}^{} \label{q:poi}\\
  &\|v^{<\kappa}\|_{\infty}^2 \le c\,\log\kappa\,\|\gb v^{<\kappa}\|_{2}^2 \label{q:kinfty}.
\end{align}
\end{lemma}

As we shall need $u\in L^\infty\cap H^1$, we consider two types of velocity.
Fixing a positive constant $U$ and denoting $k^\perp:=(-k_y,k_x)$, in section~\ref{s:sbeta} we consider (recall that $\uh_0=0$)
\begin{equation}\label{q:u}
   u(x) = \tssum_{k}\,\uh_k\ex^{\im k\cdot x} := \im U\,\tssum_{k\ne0}\,|k|^{\beta-1}k^\perp X_k\ex^{\im k\cdot x}
\end{equation}
for $\beta<-2$ (corresponding to the $|k|^{2\beta+d-1}$ energy spectrum in the physics literature).
%
In section~\ref{s:kraichnan}, we consider
\begin{equation}\label{q:ukpe}
   u(x) = \im U\,\tssum_{1\le|k|<\kpe}\,|k|^{-3}k^\perp X_k\ex^{\im k\cdot x}
\end{equation}
to handle the important marginal case $\beta=-2$ (giving the hypothesized Kraichnan's $|k|^{-3}$ energy spectrum).
This can be extended to the more realistic exponentially decaying $|\uh_k|$ for $|k|\ge\kpe$, at the cost of considerable clutter.
Further connections to 2d turbulence are discussed in \S\ref{s:kraichnan} below.

We define the complex random variable $X_k=\ex^{\im\zeta_k}$ with the random phase $\zeta_k\sim\Unif(0,2\pi)$.
This gives $|X_k|=1$ and $\EV X_k^n=0$ for any integer $n\ne0$.
We note that $X_k$ is \emph{circular\/}, i.e.\ $\EV(\ex^{\im\phi}X_k)=\EV X_k$ for any deterministic real $\phi$.
We assume that $X_k$ are uncorrelated (in the sense of complex r.v.),
\begin{equation}\label{q:EVXk}
   \EV X_j\overline{X_k}=\delta_{jk}
   \aand
   \EV X_j X_k=0
\end{equation}
subject to the constraint that $u(x)$ be real-valued, which implies $X_{-k}=\overline{X_k}$.
Therefore \eqref{q:EVXk} is taken for $j$, $k\in\Zu:=\{(k,l):l>0\}\cup\{(k,0):k>0\}$ instead of the full $\Zahl^2$.


\section{Steep energy spectrum: $\beta<-2$}\label{s:sbeta}

Noting that $k^\perp\cdot k^\perp=|k|^2$, we compute
\begin{align}
   &|\uh_k|^2 = U^2|k|^{2\beta}|X_k|^2 = U^2|k|^{2\beta}\\
   &\|u^{>\kappa}\|_2^2 = (2\pi)^2\tssum_{|k|\ge\kappa}\,|\uh_k|^2 \simeq (2\pi)^3\int_{\kappa}^\infty U^2 r^{2\beta+1} \dr \notag\\
	&\phantom{\|u^{>\kappa}\|_2^2} = \kappa^{2\beta+2}\,(2\pi)^3 U^2/(2\beta+2)
	\simeq \kappa^{2\beta+d} \|u\|_2^2 \label{q:evug}
\end{align}
where, as in \cite{jolly-dw:bht}, ``$\simeq$'' means up to $\kappa$-asymptotically smaller terms incurred by lattice approximations: $f\simeq g$ means $f=(1+(\kpb/\kappa) M)g$ with $M$ bounded; similarly for $\lesssim$.
(If a bound is needed in the above, one can replace $\kappa$ by $\kappa-1$.)
We note that even though $u$ is a random variable, thanks to $|\uh_k|$ being deterministic, these $L^2$ norms have fixed deterministic values; on the other hand, $\|u\|_p^{}$ is a random variable for $p\ne2$, which we will need to bound below.
For \eqref{q:u}, we have
\begin{align}
  \|\gb u^{>\kappa}\|_2^2 &= (2\pi)^2\tssum_{|k|\ge\kappa}\,|k|^2|\uh_k|^2 \notag\\
	&\simeq \kappa^{2\beta+4} (2\pi)^3U^2/(2\beta+4)
	= \kappa^{2\beta+2+d}\|\gb u\|_2^2.
\end{align}
Moreover, since $\beta<-2$, we have $u\in H^s\cap L^\infty$ for $s\in(d/2,-\beta-d/2)$.
We can bound $\|u\|_\infty^{}$ using Sobolev {\em including the constant\/} as
\begin{align*}
   \|u\|_\infty^2 &\le \bigl(\,\tssum_k\,|\uh_k|\bigr)^2
	& &\le \tssum_k\,|k|^{-2s}\;\tssum_k\,|k|^{2s}|\uh_k|^2\\
	&\simeq \frac{2\pi}{2s-d}\,\frac{-2\pi U^2}{2s+2\beta+d}
	& &=: f(s;\beta).
\end{align*}
The minimum of $f$ is attained for $s=-\beta/2$, giving
\begin{equation}\label{q:u8}
   \|u\|_\infty^{} \le \|\uh\|_1^{} \lesssim -2\pi U/(\beta+d)
\end{equation}
provided that $\beta+d<0$.

With these bounds on $u$, we have the following global bounds on $\tht$:

\begin{lemma}\label{t:thtg}
With $u$ given by \eqref{q:u}, $g=g^{<\kpg}$ by \eqref{q:g} and $\tht$ by \eqref{q:tht}, we have
\begin{align}
   &\|\gb\tht\|_2^{} \le \|\gbi g\|_2^{} \label{q:tht-h1}\\
   &\|\Delta\tht\|_2^{} \le c\,\|\gb u\|_2^{}\|\gb\tht\|_2^{} + \|g\|_2^{}
	\le c\,\|\gb u\|_2^{}\|\gbi g\|_2^{} + \|g\|_2^{} \label{q:tht-h2}\\
   &\|\tht\|_\infty^{} \le \|\thth\|_1^{}
	\le c\,\|\gb\tht\|_2\bigl[1 + \log^{1/2}(\|\Delta\tht\|_2^{}/\|\gb\tht\|_2^{})\bigr] \label{q:tht-l8}\\
   &\|\gb\tht\|_\infty^{} \le c\,\|\uh\|_1^{}\|\thth\|_1^{} + \|\gbi\hat g\|_1^{}. \label{q:tht-w18}
\end{align}
\end{lemma}

\noindent
We note that the $\|\gb u\|_2^{}$ in \eqref{q:tht-h2} and the $\|\uh\|_1^{}$ in \eqref{q:tht-w18} both require $\beta\le-2$, with a spectral cutoff in the marginal case.
As we do not have the lower bound for the Dirichlet quotient $\|\Delta\tht\|_2^{}/\|\gb\tht\|_2^{}$, we would have to bound $\log(\cdots)\le |\log\|\Delta\tht\|_2^{}|+|\log\|\gb\tht\|_2^{}|$ if needed.

\begin{proof}
Even though elementary for the most part, we do these explicitly as identical estimates will be used in handling \eqref{q:vfi} in~\S\ref{s:vfi}.
For \eqref{q:tht-h1}, we multiply \eqref{q:tht} by $\tht$ in $L^2$ and use the fact that $(u\cdot\gb\tht,\tht)=0$,
\begin{equation}
   \|\gb\tht\|_2^2 = (g,\tht) = (\gbi g,\agb\tht) \le \|\gbi g\|_2^{}\|\gb\tht\|_2^{}.
\end{equation}
Next, we multiply \eqref{q:tht} by $-\Delta\tht$, integrate by parts, and use the fact that $\sum_j\,(u\cdot\gb\dy_j\tht,\dy_j\tht)=0$ to get
\begin{equation}
   \|\Delta\tht\|_2^2 = -\tssum_j\,(\dy_j u\cdot\gb\tht,\dy_j\tht) - (g,\Delta\tht)
	\le \|\gb u\|_2^{}\|\gb\tht\|_4^2 + \|g\|_2^{}\|\Delta\tht\|_2^{},
\end{equation}
followed by Ladyzhenskaya, $\|\gb\tht\|_4^2\le c\,\|\gb\tht\|_2^{}\|\Delta\tht\|_2^{}$, to get \eqref{q:tht-h2}.
The inequality \eqref{q:tht-l8} is just Brezis--Gallou\"et, proved with Fourier series.
Finally, for \eqref{q:tht-w18}, we bound
\begin{equation}
   \|\gb\tht\|_\infty^{} \le \|\gb\thth\|_1^{}
	= \|\gbi\Delta\thth\|_1^{}
	= \|\gbi(u\cdot\gb\tht\hat)-\gbi\hat g\|_1
\end{equation}
and then bound the first term using $\|\gbi\gb\cdot\hat w\|_1^{}\le c\,\|\hat w\|_1^{}$ and \eqref{q:young} as
\begin{equation}
   \|\gbi\gb\cdot(u\tht\hat)\|_1^{}
	\le c\,\|(u\tht\hat)\|_1^{}
	= c\,\|\uh*\thth\|_1^{}
	\le c\,\|\uh\|_1^{}\|\thth\|_1^{}.
\end{equation}
We note that $(u\tht\hat)$ is the Fourier coefficient of $(u\tht)$, which equals the (discrete) convolution of $\uh$ and $\thth$; this type of estimate will be used again below.
\end{proof}

\noindent
Since $\|\gb u^{>\kappa}\|_2^2\sim \kappa^{2\beta+4}$, a BHT scaling for $\theta$ would mean $\EV\|\gb\tht^{>\kappa}\|_2^2\sim\kappa^{2\beta}$.
It is fairly straightforward to show a BHT-like {\em upper bound\/} for $\|\gb\tht^{>\kappa}\|_2^{}$:

\begin{prop}\label{t:bd-tht}
Let $u$ be given by \eqref{q:u} and $g=g^{<\kpg}$ by \eqref{q:g}.
Then the solution $\tht$ of \eqref{q:tht} satisfies
\begin{equation}\label{q:b-tht}
   \|\gb\tht^{>\kappa}\|_2^{} \le \kappa^\beta c(\beta)\|u\|_2^{}\|\gb\tht\|_\infty^{}
\end{equation}
for $\kappa\ge 2^\mb\max\{\kpg,\|u\|_\infty^2,c'(\beta)\|u\|_\infty^2\|u\|_2^{1/\beta}\}$ where $\mb:=\lceil-2\beta-1\rceil$.
\end{prop}

\noindent
This bound, and the identical one for $\tco$ in \eqref{q:bd-tco>} below, is suboptimal in two different ways.
First, it depends on $\|u\|_\infty^{}$ rather than $\|u\|_2^{}$ that \eqref{q:vtth-bht} suggests it should.
Second and more seriously, for its validity we need $\kappa\gtrsim c\,\|u\|^2$ rather than $\|u\|^{1/2}$ proposed by {classical theory}, at least for $\beta=-2$.
As will be seen below, this limitation re-appears in \eqref{q:bd-kpb0} and \eqref{q:bd-kpb1} in the main Theorems.

\begin{proof}
Multiplying \eqref{q:tht} by $\tht^{>\kappa}$ and using $(g,\tht^{>\kappa})=0$, we have
\begin{equation}
   \|\gb\tht^{>\kappa}\|_2^2 
	= -((u\cdot\gb\tht^{<\kappa})^{>\kappa},\tht^{>\kappa})
	\le \|\gbi(u\cdot\gb\tht^{<\kappa})^{>\kappa}\|_2^{}\|\agb\tht^{>\kappa}\|_2^{}.
\end{equation}
Now since $(u^{<\kappa/2}\cdot\gb\tht^{<\kappa/2})^{>\kappa}=0$, we can split
\begin{equation}
   (u\cdot\gb\tht^{<\kappa})^{>\kappa}
	= (u^{>\kappa/2}\cdot\gb\tht^{<\kappa})^{>\kappa} + (u^{<\kappa/2}\cdot\gb\tht_{\kappa/2,\kappa})^{>\kappa}
\end{equation}
and bound, using \eqref{q:poi}, H\"older and \eqref{q:evug},
\begin{align}
   \|\gbi&(u\cdot\gb\tht^{<\kappa})^{>\kappa}\|_2^{}
	\le \kappa^{-1}\|u^{>\kappa/2}\cdot\gb\tht^{<\kappa}\|_2^{} + \kappa^{-1}\|u^{<\kappa/2}\cdot\gb\tht_{\kappa/2,\kappa}\|_2^{} \notag\\
	&\le \kappa^{-1}\|u^{>\kappa/2}\|_2^{}\|\gb\tht^{<\kappa}\|_\infty^{} + \kappa^{-1}\|u^{<\kappa/2}\|_\infty^{}\|\gb\tht_{\kappa/2,\kappa}\|_2^{} \notag\\
	&\lesssim \kappa^\beta2^{-(\beta+1)}\|u\|_2^{}\|\gb\tht\|_\infty^{} + \kappa^{-1}\|u\|_\infty^{}\|\gb\tht^{>\kappa/2}\|_2^{}. \label{q:p2a}
\end{align}
The first term already has the correct $\kappa^\beta$ scaling, so we compute the analogous bounds for $\gb\tht^{>\kappa/2^m}$, writing $\upsilon:=\|u\|_\infty/\kappa$
\begin{align*}
   &\|\gb\tht^{>\kappa/2}\|_2^{}
	\le \kappa^\beta 2^{1-2(\beta+1)}\|u\|_2^{}\|\gb\tht\|_\infty^{} + 2\upsilon\|\gb\tht^{>\kappa/4}\|_2^{}\\
   &\|\gb\tht^{>\kappa/2^{m}}\|_2^{}
	\le \kappa^\beta2^{-(m+1)\beta-1}\|u\|_2^{}\|\gb\tht\|_\infty^{} + 2^m\upsilon\|\gb\tht^{>\kappa/2^{m+1}}\|_2^{}. \label{q:p2c}
\end{align*}
Using these repeatedly, we find
\begin{align}
   &\|\gb\tht^{>\kappa}\|_2^{} \le \kappa^\beta 2^{-\beta-1}\|u\|_2^{}\|\gb\tht\|_\infty^{}\bigl(1 + 2^{-\beta}\upsilon\bigr) + 2\upsilon^2\|\gb\tht^{\kappa/4}\|_2^{} \notag\\
	&\qquad\le \kappa^\beta 2^{-\beta-1}\|u\|_2^{}\|\gb\tht\|_\infty^{}\bigl(1 + \cdots + 2^{\mb-1-\mb\beta}\upsilon^\mb\bigr) + 2^{1+\cdots+\mb}\upsilon^{\mb+1}\|\gb\tht^{>\kappa/2^{\mb+1}}\|_2^{} \notag\\
	&\qquad\le \kappa^\beta 2^{-\beta-1}\|u\|_2^{}\|\gb\tht\|_\infty^{}\bigl(1 + \cdots + 2^{\mb-1-\mb\beta}\upsilon^\mb\bigr) + 2^{1+\cdots+\mb}\upsilon^{\mb+1}\|\gb\tht\|_2^{},
\end{align}
noting that for the penultimate step we need $(g,\tht^{>\kappa/2^\mb})=0$, i.e.\ $\kappa\ge2^\mb\kpg$.
Invoking the hypothesis $\kappa\ge2^\mb\|u\|_\infty^2$, or equivalently $\upsilon\le(2^\mb\kappa)^{-1/2}$, we can bound the bracket $(1 + \cdots + 2^{\mb-1-\mb\beta}\upsilon^\mb)$ by a constant.
Finally, demanding that the first term dominates the second,
\begin{align}
   &\kappa^\beta \|u\|_2^{}\|\gb\tht\|_\infty^{} \ge c(\beta) \kappa^{-\nu-1}\|u\|_\infty^{\nu+1}\|\gb\tht\|_2^{} \notag\\
   \Leftarrow\quad
   &\kappa^{-\beta} \ge 2^{-\mb\beta}c'(\beta)^{-\beta}\|u\|_\infty^{-2\beta}/\|u\|_2^{},
\end{align}
we recover the condition $\kappa\ge 2^\mb c'(\beta)\|u\|_\infty^2\|u\|_2^{1/\beta}$.
\end{proof}

Much of the work in the rest of this paper is to show that $\EV\|\gb\tht^{>\kappa}\|_2^2$ consists of an {\em exact expression\/} with BHT scaling plus an asymptotically smaller remainder.

\subsection{Small $U$}\label{s:ull1} 

As the current approach is different from that in \cite{jolly-dw:bht}, it is useful to consider the simpler case of small $U$ separately, as well as to motivate the analogous computation in~\S\ref{s:vtt}
We write $\tht=\tco+\thts$, which satisfy
\begin{align}
   &-\Delta\tco = g\\
   &-\Delta\thts + u\cdot\gb\thts = -u\cdot\gb\tco = u\cdot\gb\ilap g.
\end{align}
We note that since $g=g^{<\kpg}$, we also have $\tco=\tco^{<\kpg}$, so $\tht^{>\kappa}=\thts^{>\kappa}$ for $\kappa\ge\kpg$.
Then the following holds:

\begin{prop}\label{t:bht}
Assuming $\|u\|_\infty^{}\le c_*$, with $c_*<1$ defined in \eqref{q:Nb} below, and for $\kappa\ge3\kpg$, one can write $\tht=\tco+\vtt^{(1)}+\dvtt$ where
\begin{align}
   &\EV\|\gb\vtt^{(1)>\kappa}\|_2^2 \simeq \kappa^{2\beta}U^2\pi\,(2\beta)^{-1}\|\gb\tco\|_2^2\\
   &\|\gb\dvtt^{>\kappa}\|_2^{} \le \kappa^\beta\,(c + \|u\|_\infty^{})\|u\|_\infty^{}\,\|u\|_2^{}\|\gb\tco\|_\infty^{}.
\end{align}
\end{prop}

\noindent
The remainder $\gb\dvtt^{>\kappa}$ has the same $\kappa^\beta$ dependence as $\gb\vtt^{(1)>\kappa}$, but since it has a higher power of $u$, scaling as $\|u\|_\infty^{}\|u\|_2^{}\sim U^2$, it is smaller than $\bigl(\EV\|\gb\vtt^{(1)>\kappa}\|_2^2\bigr)^{1/2}$ for $\|u\|_\infty^{}$ small.

\begin{proof}
We solve for $\thts$ by iterations: put $\thts^{(0)}=0$ and
\begin{align}
   &\Delta\thts^{(n+1)} = u\cdot\gb\thts^{(n)} + u\cdot\gb\tco.
\end{align}
Now this implies
\begin{equation}\label{q:dth-b}
   \Delta\delta\thts^{(n+1)} = u\cdot\gb\delta\thts^{(n)}
\end{equation}
where $\delta\thts^{(n+1)}:=\thts^{(n+1)}-\thts^{(n)}$.
Multiplying this by $-\delta\thts^{(n+1)}$, we find
\begin{equation}
   \|\gb\delta\thts^{(n+1)}\|_2^2 = -(u\cdot\gb\delta\thts^{(n)},\delta\thts^{(n+1)})
	\le \|\gbi(u\cdot\gb\delta\thts^{(n)})\|_2^{}\,\|\gb\delta\thts^{(n+1)}\|_2^{}
\end{equation}
and thus
\begin{equation}
   \|\gb\delta\thts^{(n+1)}\|_2^{} \le \|\gbi(u\cdot\gb\delta\thts^{(n)})\|_2^{}
	\le \|u\|_\infty^{}\|\gb\delta\thts^{(n)}\|_2^{}.
\end{equation}
Therefore, if $\|u\|_\infty^{}<1$ then $\|\gb\delta\thts^{(n)}\|_2^{}\le\|u\|_\infty^{n-1}\|\gb\delta\thts^{(1)}\|_2^{}\to0$ as $n\to\infty$ and
\begin{equation}\label{q:th-1-b}\begin{aligned}
   \|\gb\thts^{(\infty)}&-\gb\thts^{(1)}\|_2^{} = \|\gb\delta\thts^{(2)}+\gb\delta\thts^{(3)}+\cdots\|_2^{}\\
	&\le \|\gb\delta\thts^{(2)}\|_2^{} + \|\gb\delta\thts^{(3)}\|_2^{} + \cdots\\
	&\le \|\gb\delta\thts^{(1)}\|_2^{}\|u\|_\infty^{}/(1 - \|u\|_\infty^{})
	= \|\gb\thts^{(1)}\|_2^{}\|u\|_\infty^{}/(1 - \|u\|_\infty^{}).
\end{aligned}\end{equation}
This tells us that $\thts=\thts^{(\infty)}$ is close to $\thts^{(1)}$ in $H^1$ if $\|u\|_\infty^{}$ is small.

Following the argument in \cite{jolly-dw:bht}, and repeated after \eqref{q:aux37} below, $\thts^{(1)}$ has BHT scaling in the sense that
\begin{equation}\label{q:th1-bht}
   \EV\|\gb\thts^{(1)>\kappa}\|_2^2 \simeq \kappa^{2\beta-2+d}\,\pi U^2 \|\gbi g\|_2^2/(2\beta+d)
\end{equation}
asymptotically for $\kappa\gg\kpg$.

In order to show that $\thts$, and thus $\tht$, also has BHT scaling asymptotically, we need to show in addition to \eqref{q:th-1-b} that $\|\gb(\thts-\thts^{(1)})^{>\kappa}\|_2^2$ is dominated by \eqref{q:th1-bht}.
We start by multiplying \eqref{q:dth-b} by $\delta\thts^{(n+1)>\kappa}$ in $L^2$,
\begin{equation}\begin{aligned}
   \|\gb\delta\thts^{(n+1)>\kappa}\|_2^2 &= (u\cdot\gb\delta\thts^{(n)},\delta\thts^{(n+1)>\kappa})
	= ((u\cdot\gb\delta\thts^{(n)})^{>\kappa},\delta\thts^{(n+1)>\kappa})\\
	&\le \|\gbi(u\cdot\gb\delta\thts^{(n)})^{>\kappa}\|_2^{}\|\agb\delta\thts^{(n+1)>\kappa}\|_2^{}.
\end{aligned}\end{equation}
To bound this, we split $u=u^{<\kappa/2}+u^{>\kappa/2}$ and $\delta\thts^{(n)<\kappa/2}+\delta\thts^{(n)>\kappa/2}$, and use the fact that $(u^{<\kappa/2}\cdot\gb\delta\thts^{(n)<\kappa/2})^{>\kappa}=0$,
\begin{align}
   \|\gbi&(u\cdot\gb\delta\thts^{(n)})^{>\kappa}\|_2^{}
	\le \kappa^{-1}\|(u\cdot\gb\delta\thts^{(n)})^{>\kappa}\|_2^{} \notag\\
	&= \kappa^{-1}\|(u^{>\kappa/2}\cdot\gb\delta\thts^{(n)} + u^{<\kappa/2}\cdot\gb\delta\thts^{(n)>\kappa/2})^{>\kappa}\|_2^{} \notag\\
	&\le \kappa^{-1}\|u^{>\kappa/2}\|_2^{}\|\gb\delta\thts^{(n)}\|_\infty^{} + \kappa^{-1}\|u^{<\kappa/2}\|_\infty^{}\|\gb\delta\thts^{(n)>\kappa/2}\|_2^{} \notag\\
	&\le c_2(\beta)\,\kappa^\beta\|u\|_2^{}\|\gb\delta\thts^{(n)}\|_\infty^{} + \kappa^{-1}\|u\|_\infty^{}\|\gb\delta\thts^{(n)>\kappa/2}\|_2^{} \label{q:bdthth>-b}
\end{align}
where we have used \eqref{q:poi} on the first line, H{\"o}lder on the third, and \eqref{q:evug} on the last.
We bound the first term as
\begin{align}
   \|\gb\delta\thts^{(n)}\|_\infty^{} &\le \|\agb^{-2}\gb(u\cdot\gb\delta\thts^{(n-1)})\|_\infty^{} & &\text{by \eqref{q:dth-b}}\notag\\
	&\le c_3\,\|u\|_\infty^{}\|\gb\delta\thts^{(n-1)}\|_\infty^{} & &\text{by \eqref{q:riesz} with $p=\infty$ and $s=1$}\notag\\
	&\le \bigl(c_3\,\|u\|_\infty^{}\bigr)^{n-1}\|\gb\thts^{(1)}\|_\infty^{} & &\text{since $\delta\vtt^{(1)}=\vtt^{(1)}$.} \label{q:u8-b}
\end{align}
We then proceed by induction, seeking to prove that there are $M_1$ and $N$ such that
\begin{equation}\label{q:dthth-mn-b}
   \|\gb\delta\thts^{(n)>\kappa}\|_2^{} \le M_1 N^{n-1} \kappa^\beta
   \qquad\text{for }\kappa\ge1 \text{ and }n\in\{1,2,3,\cdots\}.
\end{equation}
For $n=1$, we have $\delta\thts^{(1)}=\thts^{(1)}$.
We first consider $\kappa\ge2\kpg$, so $\kappa-\kpg\ge\kappa/2$,
\begin{align}
   \|\gb\thts^{(1)>\kappa}\|_2^{}
	&\le \|\gbi(u\cdot\gb\tco)^{>\kappa}\|_2^{}
	& &\le \kappa^{-1}\|(u^{>\kappa-\kpg}\cdot\gb\tco^{<\kpg})^{>\kappa}\|_2^{}\notag\\
	&\le c\,\kappa^{-1}\|u^{>\kappa-\kpg}\|_2^{}\|\gb\tco\|_\infty^{}
	& &\le c\,\kappa^{-1}\|u^{>\kappa/2}\|_2^{}\|\gb\tco\|_\infty^{}\notag\\
	&\le c_4(\beta)\,\kappa^\beta\|u\|_2^{}\|\gb\tco\|_\infty^{}, & &
\end{align}
which demands that $M_1\ge c_4(\beta)\|u\|_2^{}\|\gb\tco\|_\infty^{}$.
Then for $\kappa<2\kpg$, we bound
\begin{equation}
   \|\gb\thts^{(1)>\kappa}\|_2^{} \le \|\gb\thts^{(1)}\|_2^{} \le c_5\,\|u\|_2^{}\|\gb\tco\|_\infty^{}.
\end{equation}
Requiring that the rhs be no greater than $(2\kpg)^\beta M_1$, we find
\begin{equation}
   M_1\ge c_5\,(2\kpg)^{-\beta}\|u\|_2^{}\|\gb\tco\|_\infty^{}.
\end{equation}
Therefore \eqref{q:dthth-mn-b} holds for $n=1$ provided that $M_1\ge\max\{c_4(\beta),(2\kpg)^{-\beta}c_5\}\|u\|_2^{}\|\gb\tco\|_\infty^{}$.

Next, for $n\in\{2,3,\cdots\}$, we seek to show that each term in \eqref{q:bdthth>-b} is bounded by $\sfrac12 M_1 N^{n}\kappa^\beta$.
Using the bound
\begin{equation}\label{q:bdu-vtt1}
   \|\gb\vtt^{(1)}\|_\infty^{} = \|\agb^{-2}\gb(u\cdot\gb\tco)\|_\infty^{}
	\le c\,\|u\cdot\gb\tco\|_\infty^{}
	\qquad\text{by \eqref{q:riesz} with }s=1,
\end{equation}
we bound the first term as
\begin{align}
   c_2\|u\|_2^{}\|\gb\delta\thts^{(n)}\|_\infty^{}
	&\le c_2\|u\|_2^{}(c_3\|u\|_\infty^{})^{n-1}\|\gb\thts^{(1)}\|_\infty^{} \notag\\
	&\le c_2\|u\|_2^{}(c_3\|u\|_\infty^{})^{n-1}c_3\|u\|_\infty^{}\|\gb\tco\|_\infty^{},
\end{align}
whose rhs is bounded by $\sfrac12 N^{n-1} M_1$ if $N\ge c_3\|u\|_\infty^{}$ and $M_1\ge 2c_2c_3\|u\|_2^{}\|u\|_\infty^{}\|\gb\tco\|_\infty^{}$.
As for the second term in \eqref{q:bdthth>-b}, we bound it using \eqref{q:dthth-mn-b} as
\begin{align}
   \kappa^{-1}\|u\|_\infty\|\gb\delta\thts^{(n)>\kappa/2}\|_2^{}
	&\le \kappa^{-1}\|u\|_\infty^{} M_1 N^{n-1}(\kappa/2)^\beta \notag\\
	&= \kappa^{\beta-1} 2^{-\beta}\|u\|_\infty^{} M_1 N^{n-1}
\end{align}
whose rhs is bounded by $\sfrac12 N^n M_1 \kappa^\beta$ if $N\ge 2^{1-\beta}\|u\|_\infty^{}$.
Putting everything together, \eqref{q:dthth-mn-b} holds with
\begin{align}
   &M_1 = \max\{ c_4(\beta), (2\kpg)^{-\beta} c_5, 2c_2 c_3\|u\|_\infty \} \|u\|_2^{}\|\gb\tco\|_\infty^{},\\
   &N = \max\{ c_3, 2^{1-\beta} \} \|u\|_\infty^{}. \label{q:Nb}
\end{align}
Therefore, if $N<\sfrac12$, we have [cf.~\eqref{q:th-1-b}]
\begin{equation}\label{q:dvtt0}\begin{aligned}
   \|\gb\thts^{(\infty)>\kappa}&-\gb\thts^{(1)>\kappa}\|_2^{} = \|\gb\delta\thts^{(2)>\kappa}+\gb\delta\thts^{(3)>\kappa}+\cdots\|_2^{}\\
	&\le \|\gb\delta\thts^{(2)>\kappa}\|_2^{} + \|\gb\delta\thts^{(3)>\kappa}\|_2^{} + \cdots\\
	&\le \kappa^\beta M_1 N/(1-N) \le 2\,\kappa^\beta M_1 N.
\end{aligned}\end{equation}
This proves the proposition with $\dvtt:=\vtt^{(\infty)}-\vtt^{(1)}$.
\end{proof}


\subsection{Large $U$}

Since the above argument relies on the smallness of $\|u\|_\infty^{}$, it will need to be modified for large $U$.
For some $\kpb$ to be fixed below, we split $\tht=\tco+\vtt+\vfi$, which satisfy
\begin{align}
   &-\Delta\tco + (u\cdot\gb\tco)^{<\kpb} = g \label{q:tco}\\
   &-\Delta\vtt + (u\cdot\gb\vtt)^{>\kpb} = -(u\cdot\gb\tco)^{>\kpb} \label{q:vtt}\\
   &-\Delta\vfi + u\cdot\gb\vfi = -(u\cdot\gb\vtt)^{<\kpb}. \label{q:vfi}
\end{align}
The small $u$ case in \S\ref{s:ull1} can be regarded as a special case with $\kpb<1$ and thus $\vfi=0$, $(\cdots)^{<\kpb}=0$ and $(\cdots)^{>\kpb}=(\cdots)$.
For clarity, we shall henceforth assume that $\kpb\ge\kpg$; the general case is messier but entails no essential difficulty.

We first establish some basic facts.
Taking $\Proj^{>\kpb}$\eqref{q:tco} and using $(g^{<\kpg})^{>\kpb}=0$, we find
\begin{equation}
   \Delta\tco^{>\kpb} = 0 \Rar \tco^{>\kpb} = 0 \text{ and } \tco = \tco^{<\kpb}.
\end{equation}
This then implies $(u^{>2\kpb}\cdot\gb\tco)^{<\kpb}=(u^{>2\kpb}\cdot\gb\tco^{<\kpb})^{<\kpb}=0$, so \eqref{q:tco} is equivalent to
\begin{equation}\label{q:tco'}
   -\Delta\tco + (u^{<2\kpb}\cdot\gb\tco)^{<\kpb} = g.
\end{equation}
The fact that $\tco$ is independent of $u^{>2\kpb}$ will be important to establish the BHT scaling of $\vtt$ later.
Similarly, taking $\Proj^{<\kpb}$\eqref{q:vtt}, we have
\begin{equation}
   \Delta\vtt^{<\kpb} = 0
\end{equation}
and thus $\vtt=\vtt^{>\kpb}$.
Here we expect $\tco$ to obey the Corrsin--Obukhov scaling for (a subinterval of) $(\kpg,\kpb)$, but since we have not been able to prove it, we seek to ``isolate'' it from the rest.
As in the small $U$ case, $\vtt$ obeys the BHT scaling for large $\kappa\gg\kpb$ modulo higher-order terms, plus an additional remainder $\vfi$ which is absent when $\|u\|_\infty^{}$ is small.

The main result in this paper is the following:

\begin{theorem}\label{t:uht}
Let $u$ be given by \eqref{q:u} and $g$ by \eqref{q:g}.
We put
\begin{equation}\label{q:bd-kpb0}
   \kpb = 2^\mb\max\{\kpg, \|u\|_\infty^2,c'(\beta)\|u\|_\infty^2\|u\|_2^{1/\beta}\}
	\qquad\text{with }\mb=\lceil-2\beta-1\rceil.
\end{equation}
Then one can write $\tht=\tco^{<\kpb}+\vtt^{(1)>\kpb}+\dvtt$, where $\vtt^{(1)}$ has BHT scaling for $\kappa\ge3\kpb$,
\begin{equation}\label{q:vtth-bht}
   \EV\|\gb\vtt^{(1)>\kappa}\|_2^2 \simeq \kappa^{2\beta-2+d}\,U^2\EV\|\gb\tco\|_2^2\,4\pi^3/(-2\beta)
\end{equation}
asymptotically for $\kappa/\kpb\ll1$, while the remainder $\dvtt$ is bounded as
\begin{equation}
   \|\gb\dvtt^{>\kappa}\|_2^2 \le \kappa^\beta(M_\vtt + M_\vfi)
\end{equation}
with
\begin{align}
   &M_\vtt = c(\beta)\|u\|_2^{}\|\gb\tco\|_\infty^{}\|u\|_\infty^{}\log^2\kpb/\kpb,\\
   &M_\vfi = \kpb^\beta c(\beta)\|\uh\|_1^{}\|\gb\tco\|_\infty^{}\bigl[\|u\|_2^{} + \log^{1/2}(|\kpb^{\beta-1}\|\gb\tco\|_\infty^{})\bigr]
\end{align}
so subdominant to $\EV\|\gb\vtt^{(1)>\kappa}\|_2^2$ when $\|\uh\|_1^{}\log^2\kpb/\kpb\ll1$.
\end{theorem}

\noindent
We defer discussions to after the analogous Thm.~\ref{t:uht-ukpe}, in order to relate to classical 2d turbulence theory.


\section{Proof of Thm~\ref{t:uht}}\label{s:pf}

The proof of Theorem~\ref{t:uht} consists of three main parts.
In the first part, very similar to that of Proposition~\ref{t:bd-tht}, we show that $\|\gb\tco^{>\kpb/2}\|_2^{}\lesssim\kpb^\beta$ if $\kpb$ is large enough.
In the second part, partly resembling that of Proposition~\ref{t:bht}, we show that $\EV\|\gb\vtt^{(1)>\kappa}\|_2^2$ has a BHT scaling and that the higher-order remainder is small.
The final part, with no direct analogue in the small $u$ case, bounds $\vfi$.
The theorem follows upon taking $\delta\vtt=\vtt-\vtt^{(1)}+\vfi$.

\subsection{Bounding $\tco$} 

By multiplying \eqref{q:tco} by $\tco$ and $-\Delta\tco$ and bounding, we find the following analogues of \eqref{q:tht-h1}--\eqref{q:tht-h2},
\begin{align}
   &\|\gb\tco\|_2^{} \le \|\gbi g\|_2^{} \label{q:tco-h1}\\
   &\|\Delta\tco\|_2^{} \le c\,\|\gb u\|_2^{}\|\gb\tco\|_2^{} + \|g\|_2^{}
	\le c\,\|\gb u\|_2^{}\|\gbi g\|_2^{} + \|g\|_2^{} \label{q:tco-h2}\\
   &\|\gb\tco\|_\infty^{} \le c\,\log\kpb\,\|\Delta\tco\|_2^{}
\end{align}
where one could replace $u$ by $u^{<2\kpb}$ if desired, here and in the rest of this subsection.
Assuming that $2^{-\mb}\kpb\ge\kpg$, we then multiply \eqref{q:tco} by $\tco^{>\kappa}$ for $\kappa\in[\kpb/2,\kpb)$ to get
\begin{equation}
   \|\gb\tco^{>\kappa}\|_2^2 = -(u\cdot\gb\tco,\tco^{>\kappa})
	= -((u^{<\kappa/2}\cdot\gb\tco^{>\kappa/2}+u^{>\kappa/2}\cdot\gb\tco)^{>\kappa},\tco^{>\kappa}),
\end{equation}
giving the bound [cf.~\eqref{q:p2a}]
\begin{equation}\label{q:tco3}
   \|\gb\tco^{>\kappa}\|_2^{} \le \kappa^\beta 2^{-(\beta+1)}\|u\|_2^{}\|\gb\tco\|_\infty^{} + \kappa^{-1}\|u\|_\infty^{}\|\gb\tco^{>\kappa/2}\|_2^{}.
\end{equation}
As before, we iterate the second term to get [cf.~\eqref{q:p2c}, recalling that $\mb=\lceil-2\beta-1\rceil$ and $\upsilon=\|u\|_\infty^{}/\kappa$]
\begin{equation}\begin{aligned}
   \|\gb\tco^{>\kappa}\|_2^{} \le \kappa^\beta 2^{-\beta-1}\|u\|_2^{}\|\gb\tco\|_\infty^{}&\bigl\{1 + \cdots + 2^{\mb-1-\mb\beta}\upsilon^\mb\bigr\}\\
	&+ 2^{1+\cdots+\mb}\upsilon^\mb\|\gb\tco^{>\kappa/2^{\mb+1}}\|_2^{},
\end{aligned}\end{equation}
keeping in mind that we need $\kappa\ge\kpg$ in \eqref{q:tco3}, and thus $\kappa/2^\nu\ge\kpg$ for the last term.
Not surprisingly, the bound in Proposition~\ref{t:bd-tht} also applies to $\tco$,
\begin{equation}\label{q:bd-tco>}
   \|\gb\tco^{>\kappa}\|_2^{} \le \kappa^{\beta} c(\beta)\|u\|_2^{}\|\gb\tco\|_\infty^{}
\end{equation}
for $\kappa\in[\kpb/2,\kpb)$ provided that $\kpb\ge 2^\mb\max\{\kpg,\|u\|_\infty^2,c'(\beta)\|u\|_\infty^2\|u\|_2^{1/\beta}\}$.

\subsection{Computing $\vtt^{>\kappa}$}\label{s:vtt} 

We now show that $\vtt$ consists of a BHT-like part and a subdominant remainder.
As in Proposition~\ref{t:bht}, we solve \eqref{q:vtt} by iteration, but since $u$ is no longer small, we make use of the smallness of $u^{>\kpb}$, with $\|u\|_\infty^{}/\kpb$ in place of $\|u\|_\infty^{}$.
With $\tco$ considered given, we formally put $\vtt^{(0)}=0$ and
\begin{align}
   &\Delta\vtt^{(1)} = (u\cdot\gb\tco^{<\kpb})^{>\kpb}\\
   &\Delta\vtt^{(n+1)} = (u\cdot\gb\vtt^{(n)})^{>\kpb} + (u\cdot\gb\tco^{<\kpb})^{>\kpb}. \label{q:ivtt}
\end{align}
Unlike in the small $u$ case, we cannot obtain a BHT scaling directly since $u$ and $\tco$ are not independent, so we split $\vtt^{(1)}=\vtt^{(1)}_L+\vtt^{(1)}_H$ where
\begin{equation}
   \Delta\vtt^{(1)}_L = (u^{<2\kpb}\cdot\gb\tco^{<\kpb})^{>\kpb}
   \aand
   \Delta\vtt^{(1)}_H = (u^{>2\kpb}\cdot\gb\tco^{<\kpb})^{>\kpb}.
\end{equation}
Thanks to \eqref{q:tco'}, $\tco^{<\kpb}$ and $u^{>2\kpb}$ are independent, giving a BHT scaling for $\vtt^{(1)}_H$, as we will prove shortly.
We note the following facts:
\begin{equation}
   \vtt^{(1)}_L = \Proj_{\kpb,3\kpb}\vtt^{(1)}_L
   \aand
   \vtt^{(1)}_H = \Proj_{\kpb,\infty}\vtt^{(1)}_H.
\end{equation}
Moreover, even though $(\vtt^{(1)}_L,\vtt^{(1)}_H)_{L^2}^{}\ne0$, we do have
\begin{equation}
   \EV(\vtt^{(1)}_L,\vtt^{(1)}_H)_{L^2}^{} = 0.
\end{equation}
This can be seen by writing
\begin{equation}
   \tco(x) = \tssum_{|k|<\kpb}\,\tcoh_k \ex^{\im k\cdot x}
\end{equation}
where $\tcoh_k$ is a complex random variable that is independent of $X_j$ for $|j|\ge2\kpb$ (but {\em a priori\/} the $\tcoh_k$ are not necessarily independent mode-wise, or with $u^{<2\kpb}$); we define $\tcoh_k\equiv0$ when $|k|\ge\kpb$.
Then (writing $k\wedge j:=k_xj_y-k_yj_x$)
\begin{align*}
   \EV(\vtt^{(1)}_L,\vtt^{(1)}_H)
	= 4\pi^2U^2\tssum_{kjl}\,|k-j|^{\beta-1}&|k-l|^{\beta-1}(k\wedge j)(k\wedge l)\\
	&|j|^{-2}|l|^{-2} \EV X_{k-j} \tcoh_j \overline{X_{k-l} \tcoh_l}.
\end{align*}
But $\EV X_{k-j} \tcoh_j \overline{X_{k-l} \tcoh_l}=\EV X_{k-j} \tcoh_j \overline{\tcoh_l}\, \EV\overline{X_{k-l}}=0$ since, with $|j|$, $|l|<\kpb$, $|k-j|<2\kpb$ and $|k-l|\ge2\kpb$, $X_{k-l}$ is independent of $\tco$ and $X_{k-j}$, and $\EV X_{k-l}=0$.
We thus have
\begin{equation}
   \EV\|\vtt^{(1)}\|_2^2 = \EV\|\vtt^{(1)}_L\|_2^2 + \EV\|\vtt^{(1)}_H\|_2^2,
\end{equation}
so $\EV\|\vtt^{(1)>\kappa}\|_2^2=\EV\|\vtt^{(1)>\kappa}_H\|_2^2$ for $\kappa\ge3\kpb$.

Writing $\vttH:=\vtt^{(1)}_H$ for short, we now consider a single mode $\vttH_k$ with $|k|\ge3\kpb$,
\begin{equation}\label{q:aux37}
   \vttH_k = -U |k|^{-2} \tssum_{|j|<\kpb}\, |k-j|^{\beta-1} (k\wedge j) X_{k-j}\tcoh_j,
\end{equation}
thus
\begin{equation}
   \EV|\vttH_k|^2 = |k|^{-4}U^2\tssum_{jl}\,|k-j|^{\beta-1}|k-l|^{\beta-1}(k\wedge j)(k\wedge l)\,\EV X_{k-j}\tcoh_j \overline{X_{k-l}\tcoh_l}.
\end{equation}
Since $|k|\ge3\kpb$ and $|j|$, $|l|<\kpb$, we have $|k-j|$, $|k-l|\ge2\kpb$ and thus $X_{k-j}$ and $X_{k-l}$ are independent of the $\tco$.
Therefore $\EV X_{k-j}\tcoh_j \overline{X_{k-l}\tcoh_l}=\EV X_{k-j}\overline{X_{k-l}}\,\EV\tcoh_j\overline{\tcoh_l}=\delta_{jl}\EV|\tcoh_j|^2$, giving
\begin{equation}
   \EV|\vttH_k|^2 = |k|^{-4}U^2\tssum_j\,|k-j|^{2\beta-2}(k\wedge j)^2\EV|\tcoh_j|^2.
\end{equation}
We then sum over all $k$ with $|k|\ge\kappa$,
\begin{align*}
   \EV\|\gb\vttH^{>\kappa}\|_2^2
	&= 4\pi^2 \tssum_{|k|\ge\kappa}\,|k|^{-2}U^2\, \tssum_{|j|<\kpb}\,|k-j|^{2\beta-2}(k\wedge j)^2 \EV|\tcoh_j|^2\\
	&\simeq 4\pi^2 U^2\tssum_j\,\EV|\tcoh_j|^2\,\tssum_k\,|k|^{2\beta-4}(k\wedge j)^2\\
	&\simeq 4\pi^2 U^2\tssum_j\,|j|^2\EV|\tcoh_j|^2\,\int_0^{2\pi}\int_\kappa^\infty \sin^2\varpi_j r^{2\beta-1} \dr \;\mathrm{d}\varpi_j\\
	&= 4\pi^2 U^2 \tssum_j\,\EV|j\tcoh_j|^2 \pi\kappa^{2\beta}/(-2\beta),
\end{align*}
which gives the BHT tracer spectrum \eqref{q:vtth-bht}.
For the second line, we have approximated $|k-j|^{2\beta-2}$ by $|k|^{2\beta-2}$,
and for the third line, we have written $(k\wedge j)=:|j|\,|k|\,\sin\varpi_j$ and approximated the $k$-sum by an integral;
the errors in these approximations were shown to be $\kappa$-asymptotically subdominant in \cite{jolly-dw:bht}.

Turning to the rest of $\vtt$, we first consider global convergence.
Writing $\dvtt^{(1)}:=\vtt^{(1)}$ and $\dvtt^{(n+1)}:=\vtt^{(n+1)}-\vtt^{(n)}$ for $n=2,\cdots$, we have from \eqref{q:ivtt}
\begin{equation}\label{q:dvtt}
   \Delta\dvtt^{(n+1)} = (u\cdot\gb\dvtt^{(n)})^{>\kpb}.
\end{equation}
This then gives [cf.~\eqref{q:u8-b}]
\begin{align}
   \|\gb\dvtt^{(n+1)}\|_\infty^{}
	&= \|\agb^{-2}\gb(u\cdot\gb\dvtt^{(n)})^{>\kpb}\|_\infty^{} & &\notag\\
	&\le c_3\kpb^{-1}\log^2\kpb\,\|(u\cdot\gb\dvtt^{(n)})^{>\kpb}\|_\infty^{} & &\text{by \eqref{q:riesz-s} with $p=\infty$, $s=1$} \notag\\
	&\le c_3\kpb^{-1}\log^2\kpb\,\|u\|_\infty^{}\|\gb\dvtt^{(n)}\|_\infty^{} \notag\\
	&\le (c_3\|u\|_\infty^{}\log^2\kpb/\kpb)^n\|\gb\vtt^{(1)}\|_\infty^{} & &\text{since }\delta\vtt^{(1)}=\vtt^{(1)}\notag\\
	&\le (c_3\|u\|_\infty^{}\log^2\kpb/\kpb)^{n+1}\|\gb\tco\|_\infty^{}, \label{q:vtt-b8}
\end{align}
where for the last step we have used [cf.~\eqref{q:bdu-vtt1}], proved using \eqref{q:riesz-s} with $s=1$,
\begin{equation}
   \|\gb\vtt^{(1)}\|_\infty^{} = \|\agb^{-2}\gb(u\cdot\gb\tco)^{>\kpb}\|_\infty^{}
	\le c\kpb^{-1}\log\kpb\,\|u\cdot\gb\tco\|_\infty^{}.
\end{equation}
Therefore the iteration \eqref{q:ivtt} converges with $\vtt^{(n)}\in W^{1,\infty}$ if $c_3\|u\|_\infty^{}\log^2\kpb/\kpb<1$.

We now show that $\vtt^{>\kappa}$ is dominated by $\vtt^{(1)>\kappa}$; here the argument is very similar to that in Prop.~\ref{t:bht}, so we omit details where a computation is essentially repeated.
Noting that $\dvtt^{(n)<\kpb}=0$ by construction, we seek to prove that there are $M_1$ and $N$ such that
\begin{equation}\label{q:dhyp}
   \|\gb\dvtt^{(n)>\kappa}\|_2^{} \le M_1 N^{n-1}\kappa^\beta
	\qquad\text{for }\kappa\ge\kpb
        \text{ and }n\in\{1,2,\cdots\}.
\end{equation}
For $\dvtt^{(1)}$, we first compute, for $\kappa\ge\kpb$,
\begin{align}
   \|\gbi(u\cdot\gb\tco)^{>\kappa}\|_2^{}
	&\le \kappa^\beta c'(\beta)\|u\|_2^{}\|\gb\tco\|_\infty^{} + \kappa^{-1}\|u\|_\infty^{}\|\gb\tco^{>\kappa/2}\|_2^{} \notag\\
	&\le \kappa^\beta c(\beta)\|u\|_2^{}\|\gb\tco\|_\infty^{}\label{q:ugtco}
\end{align}
where we have assumed that $\|u\|_\infty^2/\kpb\le2^{-\mb}$ and then used \eqref{q:bd-tco>} for the last line.
This then gives
\begin{equation}
   \|\gb\vtt^{(1)>\kappa}\|_2^{} \le \kappa^\beta c(\beta)(\|u\|_2^{}\|\gb\tco\|_\infty^{}.
\end{equation}
Therefore \eqref{q:dhyp} holds for $n=1$ with $M_1\ge c(\beta)\|u\|_2^{}\|\gb\tco\|_\infty^{}$.
Next, we multiply \eqref{q:dvtt} by $-\dvtt^{(n+1)>\kappa}$ to obtain [cf.~\eqref{q:bdthth>-b}]
\begin{equation}
   \|\gb\dvtt^{(n+1)>\kappa}\|_2^{} \le c_2(\beta)\kappa^\beta\|u\|_2^{}\|\gb\dvtt^{(n)}\|_\infty^{} + \kappa^{-1}\|u\|_\infty^{}\|\gb\dvtt^{(n)>\kappa/2}\|_2^{}.
\end{equation}
Using \eqref{q:vtt-b8} and \eqref{q:dhyp}, we bound this as
\intomargin
\begin{equation}
   \|\gb\dvtt^{(n+1)>\kappa}\|_2^{}
	\le c_2\|u\|_2^{}\|\gb\tco\|_\infty^{}(c_3\|u\|_\infty^{}\log^2\kpb/\kpb)^{n}\kappa^\beta + 2^{-\beta}\kappa^{-1}\|u\|_\infty^{} M_1 N^{n-1} \kappa^\beta,
\end{equation}
which satisfies \eqref{q:dhyp} for $n+1$ provided that
\intomargin
\begin{equation}\begin{aligned}
   &M_1 = 2 \max\{c_2(\beta),c(\beta)\}\|u\|_2^{}\|\gb\tco\|_\infty^{},\\
   &N = \max\{c_3, 2^{1-\beta}\}\|u\|_\infty^{}\log^2\kpb/\kpb.
\end{aligned}\end{equation}
As in \eqref{q:dvtt0}, assuming that $N\le\sfrac12$, we have
\begin{equation}
   \|\gb(\vtt-\vtt^{(1)})^{>\kappa}\|_2^{} \le 2\kappa^\beta M_1 N =: \kappa^\beta M_\vtt
\end{equation}
which is asymptotically subdominant to $(\EV\|\gb\vtt^{(1)>\kappa}\|_2^2)^{1/2}$ for small $\|u\|_\infty\log^2\kpb/\kpb$.

\subsection{Bounding $\vfi$}\label{s:vfi} 

We note that \eqref{q:vfi} has the same form as \eqref{q:tht}, with $\vfi$ in place of $\tht$ and $\Gamma=\Gamma^{<\kpb}:=-(u\cdot\gb\vtt)^{<\kpb}$ in place of $g=g^{<\kpg}$.
We can thus apply Lemma~\ref{t:thtg} and Prop.~\ref{t:bd-tht}, {\em mutatis mutandis\/}, giving us
\begin{equation}\label{q:b-vfi}
   \|\gb\vfi^{>\kappa}\|_2^{}
	\le \kappa^\beta c(\beta)\|u\|_2^{}\|\gb\vfi\|_\infty^{}
	\le \kappa^\beta c(\beta)(\|u\|_2^{}+c)\|\gb\vfi\|_\infty^{}
\end{equation}
valid for $\kappa\ge 2^\mb\max\{\kpb,\|u\|_\infty^2,c'(\beta)\|u\|_\infty^2\|u\|_2^{1/\beta}\}$.
This already has the correct $\kappa^\beta$ scaling, so we only need to show that the coefficient $c(\beta)(\|u\|_2^{}+c)\|\gb\vfi\|_\infty^{}$ can be made smaller than the $c'(\beta)\|u\|_2^{}\|\gb\tco\|_2^{}$ in \eqref{q:vtth-bht}.

In what follows, we shall write our bounds in terms of $\tco$ as in Thm.~\ref{t:uht} rather than $g$.
We need to bound [cf.~\eqref{q:tht-w18}]
\begin{equation}\label{q:vfi4}
   \|\gb\vfi\|_\infty^{} \le c\,\|\uh\|_1^{}\|\hat\vfi\|_1^{} + \|\gbi\Gamma\|_\infty^{}.
\end{equation}
To handle the last term,
\begin{equation}
   \|\gbi\Gamma\|_\infty^{} = \|\gbi\gb\cdot(u\vtt)^{<\kpb}\|_\infty^{}
	\le c\,\|(u\vtt\hat)\|_1^{} \le c\,\|\uh*\vtth\|_1^{}
	\le c\,\|\uh\|_1^{}\|\vtth\|_1^{},
\end{equation}
recalling that $\vtt\equiv\vtt^{>\kpb}$, we use Agmon to bound [cf.~\eqref{q:aux51} for the first step]
\begin{align*}
   \|\vtth^{>\kpb}\|_1^{}
	&\le c\,\|\vtt^{>\kpb}\|_2^{1/2}\|\Delta\vtt^{>\kpb}\|_2^{1/2}
	\le c\,\kpb^{-1}\|\Delta\vtt^{>\kpb}\|_2^{}\\
	&\le c\kpb^{-1}\,\|u\|_\infty^{}\|\gb\vtt\|_2^{} + c\kpb^{-1}\|(u\cdot\gb\tco)^{>\kpb}\|_2^{}\\
	&\le c\,\kpb^{-1}\|u\|_\infty^{}\|\gbi(u\cdot\gb\tco)^{>\kpb}\|_2{} + c\kpb^{-1}\|(u\cdot\gb\tco)^{>\kpb}\|_2^{}.
\end{align*}
Noting \eqref{q:ugtco} and the analogous estimate $\|(u\cdot\gb\tco)^{>\kpb}\|_2^{}\le \kpb^{\beta+1}c(\beta)\|u\|_2^{}\|\gb\tco\|_\infty^{}$, we have
\intomargin
\begin{equation}
   \|\vtth^{>\kpb}\|_1^{}
	\le \kpb^{\beta}c(\beta)\|u\|_2^{}\|\gb\tco\|_\infty^{}(1+\|u\|_\infty/\kpb)
	\le \kpb^{\beta}c(\beta)\|u\|_2^{}\|\gb\tco\|_\infty^{}
\end{equation}
using our standing assumption that $\kpb\ge\|u\|_\infty{}$, to give
\begin{equation}\label{q:aux3vfi1}
   \|\gbi\Gamma\|_\infty^{} \le c\,\|\uh\|_1^{}\|\vtth\|_1^{}
	\le \kpb^\beta c(\beta)\,\|\uh\|_1^{}\|u\|_2^{}\|\gb\tco\|_\infty^{}.
\end{equation}
This shows that, for $u$ given, we can make the rhs of \eqref{q:vfi} as small as we like by taking $\kpb$ large.

We still need to bound the $\|\hat\vfi\|_1^{}$ in \eqref{q:vfi4}.
We first multiply \eqref{q:vfi} by $\vfi$, use \eqref{q:gbiu} and bound,
\begin{align}
   \|\gb\vfi\|_2^{} &\le \|\gbi\Gamma\|_2^{}
	\le \|(u\vtt)^{<\kpb}\|_2^{}
	\le \|u\|_\infty^{}\|\vtt^{>\kpb}\|_2{}
	\le \|\gb\vtt^{>\kpb}\|_2{}\|u\|_\infty^{}/\kpb \notag\\
	&\le \|\gbi(u\cdot\gb\tco)^{>\kpb}\|_2\|u\|_\infty^{}/\kpb
	\le \kpb^\beta c(\beta)\|u\|_2^{}\|\gb\tco\|_\infty^{}\|u\|_\infty/\kpb \notag\\
	&\le \kpb^\beta c(\beta)\|\gb\tco\|_\infty^{} \qquad\text{using }\|u\|_2^{}\|u\|_\infty^{}/\kpb\le c.
\end{align}
Next, bounding \eqref{q:vfi} directly using the above bounds and \eqref{q:ugtco},
\begin{align*}
   \|\Delta\vfi\|_2^{} &\le \|u\cdot\gb\vfi\|_2^{} + \|\Gamma\|_2^{}
	\le \|u\|_\infty^{}\|\gb\vfi\|_2^{} + \|u\cdot\gb\vtt^{>\kpb}\|_2^{}\\
	&\le \|u\|_\infty^{}(\|\gb\vfi\|_2^{} + \|\gb\vtt\|_2^{})
	\le \|u\|_\infty^{}(\|\gb\vfi\|_2^{} + \|\gbi(u\cdot\gb\tco)^{>\kpb}\|_2^{})\\
	&\le \|u\|_\infty^{}(\|\gb\vfi\|_2^{} + \kpb^\beta c(\beta)\|u\|_2^{}\|\gb\tco\|_\infty^{})\\
	&\le \kpb^\beta c(\beta)\|u\|_\infty^{}\|u\|_2^{}\|\gb\tco\|_\infty^{}(1+\|u\|_\infty^{}/\kpb)\\
	&\le \kpb^{\beta-1} c(\beta)\|\gb\tco\|_\infty^{} \qquad\text{again using }\|u\|_2^{}\|u\|_\infty^{}/\kpb\le c.
\end{align*}
Using these in [cf.\ Lemma~\ref{t:thtg}]: by \eqref{q:bgu},
\begin{align}
   \|\hat\vfi\|_1^{} &\le c\,\|\gb\vfi\|_2^{}\bigl[1 + \log^{1/2}\bigl(\|\Delta\vfi\|_2^{}/\|\gb\vfi\|_2^{}\bigr)\bigr] \notag\\
	&\le c\,\|\gb\vfi\|_2^{}\bigl[1 + \bigl(|\log\|\Delta\vfi\|_2^{}| + |\log\|\gb\vfi\|_2^{}|\bigr)^{1/2}\bigr] \notag\\
	&\le c(\beta)\kpb^\beta\|\gb\tco\|_\infty^{}\bigl[1 + \log^{1/2}(\kpb^{\beta-1}c(\beta)\|\gb\tco\|_\infty^{})\bigr]. \label{q:aux3vfi2}
\end{align}
Using \eqref{q:aux3vfi1} and \eqref{q:aux3vfi2} in \eqref{q:vfi4} and then \eqref{q:b-vfi} gives the $M_\vfi$ in the theorem.


\section{Kraichnan energy spectrum: $\beta=-2$ with cut-off}\label{s:kraichnan}

For the ``Kraichnan'' velocity \eqref{q:ukpe}, we have
\begin{align}
   &\|u\|_2^2 = (2\pi)^2U^2 \tssum_k\,|k|^{-4}
	\simeq (2\pi)^3 U^2 \int_1^{\kpe} r^{-3}\;\mathrm{d}r
	= (2\pi)^3 U^2 (1-\kpe^{-2})/2\\
   &\|u^{>\kappa}\|_2^2 \simeq (2\pi)^3 U^2 \int_{\kappa}^{\kpe} r^{-3}\;\mathrm{d}r
	= (2\pi)^3 U^2 (\kappa^{-2}-\kpe^{-2})/2
	\simeq \kappa^{-2}\|u\|_2^2\\
   &\|\gb u\|_2^2 = (2\pi)^2U^2 \tssum_k\,|k|^{-2}
	\simeq (2\pi)^3 U^2 \int_1^{\kpe} r^{-1}\;\mathrm{d}r
	= (2\pi)^3 U^2\log\kpe\\
   &\|\uh\|_1^{} = U \tssum_k\,|k|^{-2} = \|\gb u\|_2^2/(4\pi^2U).
\end{align}
Here we need to modify our ``$\simeq$'' to mean up to (absolutely bounded) terms of order $\kappa/\kpe$ as well as those of order $1/\kappa$.
We note that in what follows the only terms growing in $\kpe$ do so no faster than $\log^\alpha\kpe$.

According to the classical Batchelor--Kraichnan--Leith picture of two-dimensional turbulence \cite{kraichnan:67,leith:68,batchelor:69} (see, e.g., \cite{davidson:t,lesieur:tf,vallis:aofd} for modern synopses), in the enstrophy inertial range with enstrophy transfer rate $\eta$ the energy spectrum scales as $\eta^{2/3}|k|^{-3}$; comparing this with our $U^2|k|^{-4+1}$, we have $\eta\propto U^3$.
This enstrophy inertial range is supposed to lie between the injection scale (absent for our synthetic velocity) and the dissipation (``Kraichnan'') scale $\propto (\nu_k^3/\eta)^{1/6}$ where $\nu_k$ is the kinematic viscosity.
In the 2d Navier--Stokes equations, one can specify the forcing and viscosity independently, implying that $\|u\|_*^{}$ and $\kpe$ can be independently varied; this motivates us to prescribe $U$ and $\kpe$ in our synthetic velocity \eqref{q:ukpe}.

All the bounds in Lemma~\ref{t:thtg} also apply to this case, with $\|u\|_2^{}$, $\|\gb u\|_2^{}$ and $\|\uh\|_1^{}$ as given above (and identical proofs).
As does Prop.~\ref{t:bd-tht}: with $u$ given by \eqref{q:ukpe} and $g$ by \eqref{q:g}, we have (here $\mb=3$)
\begin{equation}
   \|\gb\tht^{>\kappa}\|_2^{} \le c\,\kappa^{-2}\|u\|_2^{}\|\gb\tht\|_\infty^{}
\end{equation}
for $\kappa\ge 2^3\max\{\kpg,\|u\|_\infty^2,c\|u\|_\infty^2\|u\|_2^{1/\beta}\}$.
The proof is again identical to that of Prop.~\ref{t:bd-tht}.

Next, we split $\tht=\tco+\vtt+\vfi$, which as before satisfy \eqref{q:tco}--\eqref{q:vfi}.
All the estimates in \S\ref{s:pf} apply with $\beta=-2$ and $\mb=3$, giving us:

\begin{theorem}\label{t:uht-ukpe}
Let $u$ be given by \eqref{q:ukpe} and $g$ by \eqref{q:g}.
We put
\begin{equation}\label{q:bd-kpb1}
   \kpb = 2^3\max\{\kpg, \|u\|_\infty^2,c\,\|u\|_\infty^2/\|u\|_2^{1/2}\}.
\end{equation}
Then one can write $\tht=\tco^{<\kpb}+\vtt^{(1)>\kpb}+\dvtt$, where $\vtt^{(1)}$ has BHT scaling for $\kappa\ge3\kpb$,
\begin{equation}\label{q:evtt-ukpe}
   \EV\|\gb\vtt^{(1)>\kappa}\|_2^2 \simeq \kappa^{-4}\,U^2\EV\|\gb\tco\|_2^2\,\pi^3/4
\end{equation}
asymptotically for $\kpb\ll\kappa\ll\kpe$, while the remainder $\dvtt$ is bounded as
\begin{equation}
   \|\gb\dvtt^{>\kappa}\|_2^{} \le \kappa^{-2}(M_\vtt + M_\vfi)
\end{equation}
with
\begin{align}
   &M_\vtt = c\,\|u\|_2^{}\|\gb\tco\|_\infty^{}\|u\|_\infty^{}\log^2\kpb/\kpb,\\
   &M_\vfi = \kpb^{-2} c\,\|\uh\|_1^{}\|\gb\tco\|_\infty^{}\bigl[\|u\|_2^{} + \log^{1/2}(\kpb^{-3}\|\gb\tco\|_\infty^{})\bigr].
\end{align}
\end{theorem}
\noindent
As in Thm.~\ref{t:uht}, here one can make the remainders $M_\vtt$ and $M_\vfi$ smaller than the BHT part $(\EV\|\gb\vtt^{(1)}\|_2^2)^{1/2}$ by taking $\kpb$ large enough.
Obviously for this Theorem to be non-empty we must have $\kpb\ll\kpe$, but as noted above, this is always possible (mathematically), as one can specify $\kpe$ independently of $U$ (and thus $\kpb$).

Now according to the classical BHT argument, the transition between the Obukhov--Corrsin and BHT regimes should take place around $\kpb_{clas}\propto\eta^{1/6}\propto U^{1/2}$.
In this regard, our $\kpb\propto U^2$ is much too large (noting that our $\kpb$ is a sufficient threshold: one is guaranteed to have BHT spectrum for $\kappa\ge3\kpb$, though of course one might have it for smaller $\kappa$).
A more careful analysis might give a rigorous $\kpb\propto U^{1+\eps}$, but we do not see a way to bring the power on $U$ below~1.



\appendix
\section{Proofs of basic inequalities}\label{s:ab}

For \eqref{q:gbiu}--\eqref{q:bgu}, we define $v:=\gbi\gb\cdot u$ and note that $\vh_k=\im(k/|k|)\cdot\uh_k$.
Now for each $k$,
\[
   |\vh_k|^2 = |(k/|k|)\cdot\uh_k|^2 \le |(k/|k|)|^2|\uh_k|^2 = |\uh_k|^2
\]
which gives us \eqref{q:gbiu},
\begin{equation}
   \|v\|_2^{} = \|\vh\|_2^{} \le \|\uh\|_2^{} = \|u\|_2^{}.
\end{equation}
Similarly, we have $\|\agb^s v\|_2^{}\le c\,\|\agb^s u\|_2^{}$ for any $s\in\Real$ whenever the expressions are defined.

As noted, the next two are simply Agmon and Brezis--Gallou\"et, proved using Fourier series.
For both, we start by writing for some $\kappa\ge1$ to be determined below (recall that all wavenumber sums exclude $k=0$),
\begin{equation}
   \|v\|_\infty^{} \le \|\vh\|_1^{} = \|\vh^{<\kappa}\|_1^{} + \|\vh^{>\kappa}\|_1^{}.
\end{equation}
For \eqref{q:agmon}, we bound
\begin{equation}\label{q:aux51}
   \|\vh^{<\kappa}\|_1^{} + \|\vh^{>\kappa}\|_1^{}
	\le c\kappa\,\|\vh^{<\kappa}\|_2^{} + c\kappa^{-1}\||\cdot|^2\vh^{>\kappa}\|_2^{}
\end{equation}
and then choose $\kappa=\|\Delta v\|_2^{1/2}/\|v\|_2^{1/2}=\||\cdot|^2\vh\|_2^{1/2}/\|\vh\|_2^{1/2}$.
For \eqref{q:bgu}, we assume that $\|\gb v\|_2^{}=1$ (otherwise replace $v$ by $v/\|\gb v\|_2^{}$), bound
\begin{equation}
   \|\vh^{<\kappa}\|_1^{} + \|\vh^{>\kappa}\|_1^{}
	\le c(\log\kappa)^{1/2}\||\cdot|\vh^{<\kappa}\|_2^{} + c\kappa^{-1}\||\cdot|^2\vh^{>\kappa}\|_2^{}
\end{equation}
and choose $\kappa=\|\Delta v\|_2^{}$.

In Lemma~\ref{t:P}, the first three are obvious from the Fourier representation.
As for \eqref{q:kinfty}, we compute using Cauchy--Schwarz,
\begin{align}
   \|v^{<\kappa}\|_{\infty}^{} &\le \tssum_{|k|<\kappa}'\,|\vh_k|
	\le \tssum_{|k|<\kappa}'\,|\vh_k|\,|k|\,|k|^{-1} \notag\\
	&\le \Bigl(\tssum_{|k|<\kappa}'\,|\vh_k|^2|k|^2\Bigr)^{1/2}\Bigl(\tssum_{|k|<\kappa}'\,|k|^{-2}\Bigr)^{1/2},
\end{align}
and bound the last sum by an integral,
\begin{equation}
   \tssum_{|k|<\kappa}\,|k|^{-2} \le 2\pi\,\int_1^\kappa r^{-1}\dr = 2\pi\log\kappa.
\end{equation}


\section{Proof Lemma~\ref{t:riesz}}\label{s:riesz}

\newcommand{\That}{\hat T}
\newcommand{\Rhat}{\hat R}
\newcommand{\what}{\hat w}
\newcommand{\bgt}{\triangleright}
\newcommand{\Dtil}{\tilde D}

C.~Demeter sketched the proof of \eqref{q:riesz}, essentially in complete form.
We are responsible for the details and the (slight) generalisation to \eqref{q:logb} and \eqref{q:riesz-s}.

\begin{proof}
It is instructive to start with the one-dimensional version (a standard exercise in harmonic analysis).
Let
\begin{equation}
   w(x) = \tssum_{|k|\ge\kappa}\,\what_k \ex^{\im kx}
\end{equation}
and consider [abuse of notation: $T$ denotes both the integral operator and its kernel]
\begin{equation}
   Tw(x) = (T*w)(x)
   \quad\text{ where }\quad
   T(x) = \tssum_{|k|\ge\kappa}\,|k|^{-s}\ex^{\im kx}.
\end{equation}
Using the Dirichlet kernel
\begin{equation}
   D_N(x) := \tssum_{|k|\le N}\,\ex^{\im kx},
\end{equation}
we compute
\begin{align*}
   T(x)
	&= \tssum_{k\ge\kappa}\,k^{-s}\bigl(\ex^{\im kx}+\ex^{-\im kx}\bigr)
	= \tssum_{k\ge\kappa}\,k^{-s}[D_k(x)-D_{k-1}(x)]\\
	&= \tssum_{k\ge\kappa}\,D_k(x)\bigl(k^{-s}-(k+1)^{-s}\bigr) - D_{\kappa-1}(x)\kappa^{-s}.
\end{align*}
Now by Taylor's theorem, $k^{-s}-(k+1)^{-s}=s\,(k+\xi)^{-s-1}\le s\,k^{-s-1}$ with $\xi\in[0,1]$, so taking $\|\cdot\|_1^{}$ and using the triangle inequality,
\begin{equation}
   \|T\|_1^{} \le \tssum_{k\ge\kappa} \|D_k\|_1^{} s\,k^{-s-1} + \|D_{\kappa-1}\|_1^{}\kappa^{-s}.
\end{equation}
Using the well-known fact that
\begin{equation}
   \|D_k\|_1^{} \le c\,\log k
\end{equation}
and bounding the sum by the corresponding integral, we arrive at
\begin{equation}
   \|T\|_1^{} \le c\,\kappa^{-s}\log\kappa + s\int_{\kappa-1}^\infty k^{-s-1}\log(k+1)\;\mathrm{d}k \le c\,\kappa^{-s} L_1(\kappa,s).
\end{equation}
Young's (convolution) inequality, then implies, for all $p\in[1,\infty]$,
\begin{equation}
   \|T*w\|_p^{} \le \|T\|_1^{}\|w\|_p^{}
	\le c\,\kappa^{-s}L_1(\kappa,s)\|w\|_p^{}.
\end{equation}
An identical bound obtains for $\hat R=\im k\,|k|^{-s-1}$ by summing $k_1>0$ and $k_1<0$ separately, and using $\Dtil_N$ defined in \eqref{q:Dtil} below instead of $D_N$; we shall do this explicitly in the 2d case.

For the 2d case, we will need to separate $x=(x_1,x_2)$ and $k=(k_1,k_2)$; we write $|k|:=(k_1^2+k_2^2)^{1/2}$.
For $w\in L^2(\Dom)$ and $\kappa\ge2$, we define its high-mode ``square'' projection,
\begin{equation}
   w^{\bgt\kappa}(x) := \tssum_{|k_1|\ge\kappa,\,|k_2|\ge\kappa}\,\hat w_k^{}\ex^{\im k\cdot x}.
\end{equation}
We note that $(w^{>\kappa})^{\bgt(\kappa/\sqrt{}\,2)}=w^{>\kappa}$, so for \eqref{q:riesz-s} it suffices to prove
\begin{equation}\label{q:Rbgt}
   \|\agb^{-s-1}\gb w^{\bgt(\kappa/\sqrt{}\,2)}\|_p^{} \le c'(s)s^{-1}L_2(2\kappa,s)\|w^{\bgt(\kappa/\sqrt{}\,2)}\|_p^{}
\end{equation}
and similarly for the other inequalities.

We consider the operator $R_1=\agb^{-s-1}\dy_1$.
Instead of writing its symbol $\im k_1|k|^{-s-1}$ explicitly, we write
\begin{equation}
   f_{k_1 k_2} = f(k_1,k_2) := k_1/(k_1^2+k_2^2)^{(s+1)/2}
\end{equation}
and note the symmetry properties
\begin{equation}\label{q:Rfkk}
   f(-k_1,k_2)=-f(k_1,k_2)
   \aand
   f(k_1,-k_2)=f(k_1,k_2),
\end{equation}
as well as the bounds, regarding $k_1$ and $k_2$ as continuous variables in $\Real^2$,
\begin{equation}\label{q:bdfkk}
   |f(k_1,k_2)| \le |k|^{-s},\quad
   |\dy_i f(k_1,k_2)| \le c(s)|k|^{-s-1} \text{ and }
   |\dy_{ij} f(k_1,k_2)| \le c(s)|k|^{-s-2}
\end{equation}
for $i,\,j\in\{1,2\}$.

For $t\in[0,2\pi)$, we put
\begin{equation}\label{q:Dtil}
   \Dtil_N(t) := \tssum_{k=1}^N\,\ex^{\im kt}
	= \ex^{\im t}\frac{\ex^{\im Nt}-1}{\ex^{\im t}-1}.
\end{equation}
Even though $\Dtil_N(t)$ is not quite the usual Dirichlet kernel, we do have
\begin{equation}
  |\Dtil_N(t)| = \left|\frac{\ex^{\im Nt}-1}{\ex^{\im t}-1}\right|
	= \frac{(1-\cos Nt)^{1/2}}{(1-\cos t)^{1/2}}
	= \frac{|\sin Nt/2|}{|\sin t/2|},
\end{equation}
so it can be bounded like the Dirichlet kernel,
\begin{equation}
   \|\Dtil_N\|_{L^1(\Dom)}^{} \le c\,\log N \qquad\text{for }N\ge2.
\end{equation}
For convenience (when appearing in sums), we define $\Dtil_0=0$.
We partition our summation domain as (see figure)
\begin{align*}
   \{(k_1,k_2)&:|k_1|\ge\kappa,\,|k_2|\ge\kappa\} - \{(0,k_2):|k_2|\ge\kappa\}\\
	&= K_{a+} \cup K_{b+} \cup K_{c+} \cup K_{a-} \cup K_{b-} \cup K_{c-}
\end{align*}
\begin{figure}[ht]\label{Kfig}
\psfrag{k}{$\kappa$}
\psfrag{k1}{$k_1$}
\psfrag{k2}{$k_2$}
\psfrag{a}{$K_{a+}$}
\psfrag{b}{$K_{b+}$}
\psfrag{c}{$K_{c+}$}
\centerline{\includegraphics[scale=.5] {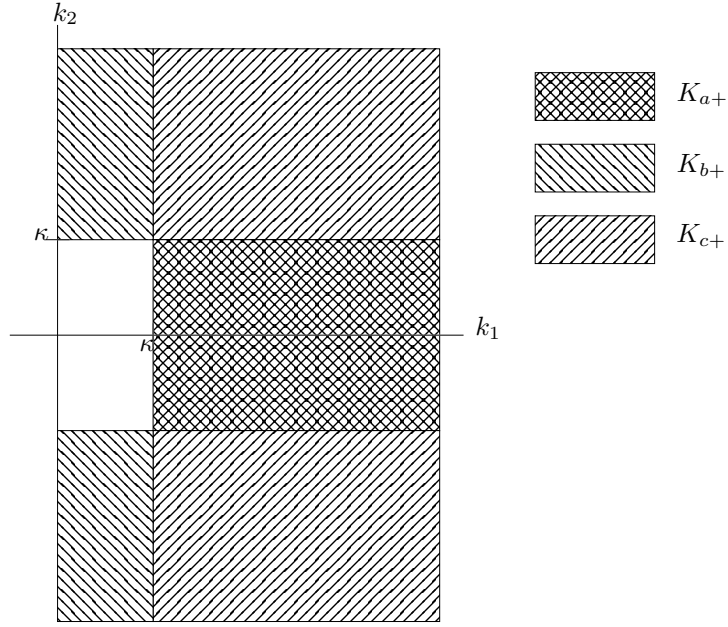}}
\caption{Partition of summation domain}
\end{figure}
where
\begin{align*}
   &K_{a+} = \{(k_1,k_2):k_1\ge\kappa,\,|k_2|\le\kappa\},\\
   &K_{b+} = \{(k_1,k_2):0<k_1\le\kappa,\,|k_2|\ge\kappa\},\\
   &K_{c+} = \{(k_1,k_2):k_1>\kappa,\,|k_2|>\kappa\},
\end{align*}
plus their analogues for $k_1<0$.
Four points $(\pm\kappa,\pm\kappa)$ are counted twice, so they are subtracted off in \eqref{q:Rn}.

We have
\begin{equation}\label{q:Rn}\begin{aligned}
   R_1(x) &= \tssum_{|k_1|,\,|k_2|\ge\kappa}\,\ex^{\im(k_1x_1+k_2x_2)}f(k_1,k_2)\\
	&= R_{a+}(x) + \cdots + R_{c-}(x) - \ex^{\im\kappa(x_1+x_2)}f_{\kappa\kappa} - \cdots - \ex^{-\im\kappa(x_1+x_2)}f_{-\kappa-\kappa}
\end{aligned}\end{equation}
where $R_{a+}(x)$ denotes the corresponding sum over $K_{a+}$, etc.
Taking $L^1$ norm, this gives us (noting that $|f_{\pm\kappa,\pm\kappa}|=|f_{\kappa\kappa}|$)
\begin{equation}
   \|R_1\|_1^{} \le \|R_{a+}\|_1^{} + \cdots \|R_{c-}\|_1^{} + 4c\,|f_{\kappa\kappa}|.
\end{equation}

Bounding each term, we start by computing
\begin{equation}
   R_{c+}(x) = \tssum_{k_1>\kappa}\,\ex^{\im k_1x_1}\tssum_{|k_2|>\kappa}\,\ex^{\im k_2x_2} f_{k_1k_2}
	=: \tssum_{k_1=\kappa+1}^\infty\,\ex^{\im k_1x_1} S_{c+,k_1}(x_2).
\end{equation}
Now using $f_{k_1,-k_2}=f_{k_1k_2}$ and $\ex^{\im ky}+\ex^{-\im ky}=D_k(y)-D_{k-1}(y)$, we have
\begin{align*}
   S_{c+,k_1}(x_2) &:= \tssum_{k_2=\kappa+1}^\infty\,(\ex^{\im k_2x_2} + \ex^{-\im k_2x_2}) f_{k_1k_2}\\
	&= \tssum_{k_2=\kappa+1}^\infty\,[D_{k_2}(x_2)-D_{k_2-1}(x_2)] f_{k_1k_2}\\
	&= \tssum_{k_2=\kappa+1}^\infty\,D_{k_2}(x_2)[f_{k_1k_2}-f_{k_1,k_2+1}] - D_{\kappa}(x_2)f_{k_1,\kappa+1}.
\end{align*}
Similarly, using $\ex^{\im kx}=\Dtil_{k}(x)-\Dtil_{k-1}(x)$, we compute (abbreviating $S_{c+,k}$ to $S_k$)
\begin{align*}
   R_{c+}(x) &= \tssum_{k_1=\kappa+1}^\infty\, \ex^{\im k_1x_1} S_{k_1}(x_2)\\
	&= \tssum_{k_1=\kappa+1}^\infty [\Dtil_{k_1}(x_1)-\Dtil_{k_1-1}(x_1)] S_{k_1}(x_2)\\
	&= \tssum_{k_1=\kappa+1}^\infty \Dtil_{k_1}(x_1) [S_{k_1}(x_2)-S_{k_1+1}(x_2)] - \Dtil_\kappa(x_1) S_{\kappa+1}(x_2)\\
	&= {\tssuml_{k_1>\kappa}}\Dtil_{k_1}(x_1)\Bigl\{ D_\kappa(x_2)[f_{k_1+1,\kappa+1}-f_{k_1,\kappa+1}]\\
	&\qquad\qquad+{\tssuml_{k_2>\kappa}}D_{k_2}(x_2)[f_{k_1k_2}-f_{k_1k_2+1}-f_{k_1+1,k_2}+f_{k_1+1,k_2+1}]\Bigr\}\\
	&\qquad+\Dtil_\kappa(x_1) \Bigl\{ D_\kappa(x_2) f_{\kappa+1,\kappa+1} - {\tssuml_{k_2>\kappa}} D_{k_2}(x_2)[f_{\kappa+1,k_2}-f_{\kappa+1,k_2+1}] \Bigr\}.
\end{align*}
Using Taylor's theorem, we can write the $f$-differences as
\begin{align*}
   &f_{k_1+1,\kappa+1}-f_{k_1,\kappa+1} = \dy_1 f(k_1+\xi,\kappa+1),\\
   &f_{\kappa+1,k_2}-f_{\kappa+1,k_2+1} = -\dy_2 f(\kappa+1,k_2+\xi),\\
   &f_{k_1k_2}-f_{k_1k_2+1}-f_{k_1+1,k_2}+f_{k_1+1,k_2+1} = \dy_{22}f(k_1+\xi,k_2+\xi'') + \dy_{12}f(k_1+\xi,k_2+\xi'')
\end{align*}
for some $\xi$, $\xi''\in[0,1]$ which may be different each time they appear.
Recalling that $\|D_k\|_1^{}\le c\,\log k$ and $\|\Dtil_k\|_1^{}\le c\log k$, we bound $\|R_{c+}\|_1^{}$ term-by-term as follows,
\begin{align}
   &\|\Dtil_\kappa(\cdot) D_\kappa(\cdot) f_{\kappa+1,\kappa+1}\|_1
	\le \|\Dtil_\kappa\|_1^{}\,\|D_\kappa\|_1^{}\,|f_{\kappa+1,\kappa+1}|
	\le c \kappa^{-s}\log^2\kappa\\
   &\bigl\|\tssum_{k_1>\kappa}\,\Dtil_{k_1} D_\kappa [f_{k_1+1,\kappa+1}-f_{k_1,\kappa+1}]\bigr\|_1^{}
	\le \tssum_{k_1>\kappa}\,\|\Dtil_{k_1}\|_1^{} \|D_\kappa\|_1^{} |\dy_1 f(k_1+\xi,\kappa+1)| \notag\\
	&\qquad\le c\log\kappa\,\tssum_{k_1>\kappa}\,k_1^{-s-1}\log k_1
	\le c\,\log\kappa\, L_1(\kappa,s)/(s\kappa^s),\\
   &\bigl\|\tssum_{k_2>\kappa}\,\Dtil_\kappa D_{k_2} [f_{\kappa+1,k_2}-f_{\kappa+1,k_2+1}]\bigr\|_1^{}
	\le \tssum_{k_2>\kappa}\,\|\Dtil_{\kappa}\|_1^{} \|D_{k_2}\|_1^{} |\dy_2 f(\kappa+1,k_2+\xi)| \notag\\
	&\qquad\le c\log\kappa\,\tssum_{k_2>\kappa}\,k_2^{-s-1}\log k_2
	\le c\,\log\kappa\,L_1(\kappa,2)/(s\kappa^s).
\end{align}
Here we have used the facts that, for $k_1>0$ and $k_2>0$, our bounds \eqref{q:bdfkk} are all monotone decreasing in $|k|$, and that $k^{-s-1}\log k$ is monotone decreasing for $k\ge\ex$ (thus requiring $\kappa\ge4\ge\sqrt2\,\ex$ in the lemma), to bound the sum as
\begin{equation}
   \tssum_{k=\kappa+1}^\infty\,k^{-s-1}\log k \le \int_{\kappa}^\infty k^{-s-1}\log k\;\mathrm{d}k = L_1(\kappa,s)/(s\kappa^s).
\end{equation}
We bound the double sum as
\begin{align*}
   &\bigl\|\tssum_{k_1,\,k_2=\kappa+1}^\infty\,\Dtil_{k_1}\,D_{k_2}\,[f_{k_1k_2}-f_{k_1k_2+1}-f_{k_1+1,k_2}+f_{k_1+1,k_2+1}]\bigr\|_1^{}\\
	&\qquad\le 2\,\tssum_{k_1,\,k_2=\kappa+1}^\infty\,\|\Dtil_{k_1}\|_1^{}\,\|D_{k_2}\|_1^{}\,|\dy_{i2}f(k_1+\xi,k_2+\xi'')|\\
	&\qquad\le c(s)\tssum_{k_1,\,k_2=\kappa+1}^\infty\,\log k_1 \log k_2\,|(k_1,k_2)|^{-s-2}\\
	&\qquad\le c(s) \int_\kappa^\infty\!\!\int_\kappa^\infty\,|(k_1,k_2)|^{-s-2}\log k_1\,\log k_2\;\mathrm{d}k_1\,\mathrm{d}k_2,
\end{align*}
followed by a change-of-variables to $n_1=\sfrac12(k_1+k_2)$, $n_2=\sfrac12(k_2-k_1)$,
\begin{align}
	&\qquad\le c(s) \int_\kappa^\infty\!\!\int_{\kappa-n_1}^{n_1-\kappa}\,|(n_1,n_2)|^{-s-2}\log(n_1+n_2)\,\log(n_1-n_2)\;\mathrm{d}n_1\,\mathrm{d}n_2\notag\\
	&\qquad\le c(s) \int_\kappa^\infty 2n_1\,|(n_1,n_1)|^{-s-2}\log(2n_1)\,\log(2n_1)\;\mathrm{d}n_1\notag\\
	&\qquad\le c(s) \int_\kappa^\infty n_1^{-s-1}\log^2(2n_1)\;\mathrm{d}n_1
	\le c(s)\,L_2(\kappa,s)/(s\kappa^s).
\end{align}
These give us, noting that $\log\kappa\,L_1(\kappa,s)\le L_2(\kappa,s)$,
\begin{equation}
   \|R_{c+}\|_1^{} \le c(s)s^{-1}\kappa^{-s}L_2(\kappa,s).
\end{equation}

Next, we compute
\begin{equation}
   R_{b+}(x) = \tssum_{|k_2|\ge\kappa}\,\ex^{\im k_2x_2} \tssum_{k_1=1}^\kappa\,\ex^{\im k_1x_1} f_{k_1k_2}
	=: \tssum_{|k_2|\ge\kappa}\,\ex^{\im k_2x_2} S_{b+,k_2}(x_1).
\end{equation}
Now using the definition $\Dtil_0=0$,
\begin{align*}
   S_{b+,k_2}(x_1) &= \tssum_{k_1=1}^\kappa [\Dtil_{k_1}(x_1)-\Dtil_{k_1-1}(x_1)] f_{k_1k_2}\\
	&= \tssum_{k_1=1}^{\kappa-1}\,\Dtil_{k_1}(x_1)[f_{k_1k_2}-f_{k_1+1,k_2}] + \Dtil_\kappa(x_1) f_{\kappa k_2},
\end{align*}
so, again writing $S_{k_2}$ for $S_{b+,k_2}$,
\begin{align*}
   R_{b+}(x) &= \tssum_{k_2=\kappa}^\infty\,[D_{k_2}(x_2)-D_{k_2-1}(x_2)] S_{k_2}(x_1)\\
	&= \tssum_{k_2=\kappa}^\infty\,D_{k_2}(x_2) [S_{k_2}(x_1)-S_{k_2+1}(x_1)] - D_{\kappa-1}(x_2)S_\kappa(x_1)\\
	&= {\tssum_{k_2=\kappa}^\infty}\,D_{k_2}(x_2)\Bigl\{ \Dtil_\kappa(x_1) [f_{\kappa k_2}-f_{\kappa k_2+1}]\\
	&\qquad\qquad + {\tssuml_{k_1=1}^{\kappa-1}}\,\Dtil_{k_1}(x_1)[f_{k_1k_2}-f_{k_1+1,k_2}-f_{k_1k_2+1}+f_{k_1+1,k_2+1}]\Bigr\}\\
	&\qquad- D_{\kappa-1}(x_2)\Bigl\{ {\tssum_{k_1=1}^{\kappa-1}}\,\Dtil_{k_1}(x_1)[f_{k_1\kappa}-f_{k_1+1,\kappa}] + \Dtil_\kappa(x_1) f_{\kappa\kappa} \Bigr\}.
\end{align*}
We can bound $\|R_{b+}\|_1^{}$ in the same way to that of $R_{c+}$ above, so we only sketch the most difficult term: since $|f_{k_1k_2}-f_{k_1+1,k_2}-f_{k_1k_2+1}+f_{k_1+1,k_2+1}|\le c\,|(k_1,k_2)|^{-s-2}$, the double sum is bounded in $L^1$ by
\begin{align*}
   {\tssuml_{k_2=\kappa}^\infty}\,&\log k_2 {\tssuml_{k_1=1}^{\kappa-1}}\,\log k_1 |(k_1,k_2)|^{-s-2}
	\le {\tssuml_{k_2=\kappa}^\infty}\,\log k_2{\tssuml_{k_1=1}^{\kappa-1}}\,\log(\kappa-1)|(1,k_2)|^{-s-2}\\
	&\le (\kappa-1)\log(\kappa-1)\tssum_{k_2=\kappa}^\infty\,\log k_2 |k_2|^{-s-2}
	\le c\,L_1(\kappa-1,s+1)\log(\kappa-1)/(s(\kappa-1)^s).
\end{align*}
The single sums being easier to bound, this gives us
\begin{equation}
    \|R_{b+}\|_1^{} \le c(s)\log(\kappa-1)L_1(\kappa,s+1)/(s(\kappa-1)^s),
\end{equation}
which is subdominant to $\|R_{c+1}\|_1^{}$.

Finally, we compute
\begin{equation}
   R_{a+}(x) = \tssum_{k_1=\kappa}^\infty\,\ex^{\im k_1x_1} \tssum_{k_2=-\kappa}^\kappa \ex^{\im k_2x_2} f_{k_1k_2}
	=: \tssum_{k_1=\kappa}^\infty\,\ex^{\im k_1x_1} S_{a+,k_1}(x_2).
\end{equation}
As before, using $f_{k_1,-k_2}=f_{k_1k_2}$,
\begin{align*}
   S_{k_1}(x_2) &= \tssum_{k_2=1}^\kappa\,(\ex^{\im k_2x_2}+\ex^{-\im k_2x_2}) f_{k_1k_2} + f_{k_10}\\
	&= \tssum_{k_2=1}^\kappa\,[D_{k_2}(x_2)-D_{k_2-1}(x_2)] f_{k_1k_2} + f_{k_10}\\
	&= {\tssuml_{k_2=1}^{\kappa-1}} D_{k_2}(x_2)[f_{k_1k_2}-f_{k_1k_2+1}] + D_\kappa(x_2) f_{k_1\kappa} - D_0(x_2)f_{k_11} + f_{k_10},
\end{align*}
giving us
\begin{align*}
   R_{a+}(x) &= \tssum_{k_1=\kappa}^\infty\,[\Dtil_{k_1}(x_1)-\Dtil_{k_1-1}(x_1)] S_{k_1}(x_2)\\
	&= \tssum_{k_1=\kappa}^\infty\,\Dtil_{k_1}(x_1)[S_{k_1}(x_2) - S_{k_1+1}(x_2)] - \Dtil_{\kappa-1}(x_1)S_\kappa(x_2)\\
	&= {\tssuml_{k_1=\kappa}^\infty}\,\Dtil_{k_1}(x_1)\bigl\{ D_\kappa(x_2)[f_{k_1\kappa}-f_{k_1+1,\kappa}+f_{k_10}-f_{k_1+1,0}-f_{k_11}+f_{k_1+1,1}]\\
	&\qquad\qquad+ \tssum_{k_2=1}^{\kappa-1}\,D_{k_2}(x_2)[f_{k_1k_2}-f_{k_1k_2+1}-f_{k_1+1,k_2}+f_{k_1+1,k_2+1}]\\
	&\qquad+ \Dtil_{\kappa-1}(x_1)\Bigl\{ {\tssuml_{k_2=1}^{\kappa-1}}\,D_{k_2}(x_2)[f_{\kappa k_2}-f_{\kappa k_2+1}] + D_\kappa(x_2)f_{\kappa\kappa} - f_{\kappa1} + f_{\kappa0}\Bigr\}.
\end{align*}
The double sum can be bounded in $L^1$ as we did with $R_{b+}$ above, while in the first sum
\begin{align*}
   |f_{k_1\kappa}&-f_{k_1+1,\kappa}+f_{k_10}-f_{k_1+1,0}-f_{k_11}+f_{k_1+1,1}|\\
	&\le |\dy_1 f(k_1+\xi,\kappa)| + c(s)\,|\dy_{ij}f(k_1+\xi,\xi'')|
	\le |k_1|^{-s-1} + c(s)|k_1|^{-s-2}
\end{align*}
so can be bounded as before.
Not surprisingly, $R_{a+}$ is bounded like $R_{b+}$,
\begin{equation}
    \|R_{a+}\|_1^{} \le c(s)s^{-1}(\kappa-1)^{-s}\log(\kappa-1)L_1(\kappa-1,s+1).
\end{equation}
Obviously identical bounds hold for $R_{a-}$, $R_{b-}$ and $R_{c-}$, so noting that $\log\kappa\,L_1(\kappa,s+1)$ is dominated by $L_2(\kappa,s)$, we have
\begin{equation}
   \|R_1\|_1^{} \le c(s)/(s\kappa^s) L_2(\kappa,2).
\end{equation}
The same bound holds for $R_2=\agb^{-s-1}\dy_y$ and, modulo a factor of 2, for $R=\agb^{-s-1}\gb$.
The operator bound \eqref{q:Rbgt} then follows via Young's inequality,
\begin{equation}
   \|R w^{\bgt\kappa'}\|_p^{} = \|R*w^{\bgt\kappa'}\|_p^{}
	\le \|R\|_1^{}\|w^{\bgt\kappa'}\|_p^{}.
\end{equation}

For \eqref{q:riesz}, we could in principle add to \eqref{q:riesz-s} the contributions of the low modes with $|k_1|$, $|k_2|<5$.
For a better bound, however, we write $R_1=R_++R_-$ where
\begin{equation}
   R_+(x) = \tssum_{k_1=1}^\infty\,\ex^{\im k_1x_1}\,\tssum_{k_2=-\infty}^\infty\,\ex^{\im k_2x_2}f_{k_1k_2}
	=: \tssum_{k_1=1}^\infty\,\ex^{\im k_1x_1}S_{k_1}(x_2).
\end{equation}
Now as before
\begin{align*}
   S_{k_1}(x_2) &= f_{k_10} + \tssum_{k_2=1}^\infty\,(\ex^{\im k_2x_2}+\ex^{-\im k_2x_2})f_{k_1k_2}\\
	&= f_{k_10} + \tssum_{k_2=1}^\infty\,[D_{k_2}(x_2)-D_{k_2-1}(x_2)] f_{k_1k_2}\\
	&= f_{k_10} - f_{k_11} + \tssum_{k_2=1}^\infty\,D_{k_2}(x_2) [f_{k_1k_2}-f_{k_1,k_2+1}]
\end{align*}
and
\begin{align*}
   R_+(x) &= \tssum_{k_1=1}^\infty\,[\Dtil_{k_1}(x_1)-\Dtil_{k_1-1}(x_1)] S_{k_1}(x_2)\\
	&= \tssum_{k_1=1}^\infty\,\Dtil_{k_1}(x_1)[S_{k_1}(x_2)-S_{k_1+1}(x_2)]\\
	&= {\tssum_{k_1=1}^\infty}\,\Dtil_{k_1}(x_1)\bigl\{ f_{k_10} - f_{k_11} - f_{k_1+1,0} + f_{k_1+1,1}\\
	&\qquad\qquad+ {\tssum_{k_2=1}^\infty}\,D_{k_2}(x_2)[f_{k_1k_2} - f_{k_1k_2+1} - f_{k_1+1,k_2} + f_{k_1+1,k_2+1}]\bigr\}.
\end{align*}
The double sum can be bounded as, up to a constant $c(s)$,
\begin{align*}
   \tssum_{k_1=1}^\infty\,&\log k_1\,\tssum_{k_2=1}^\infty\,|(k_1,k_2)|^{-s-2}\log k_2
	\le \tssum_{k_1=1}^\infty\tssum_{k_2=1}^\infty\,|k|^{-s-2}\log^2|k|\\
	&\le \frac\pi2 \int_1^\infty k^{-s-1}\log^2 k\;\mathrm{d}k
	= c/s^3.
\end{align*}
The single sum is subdominant, and the same bound applies to $R_-$, giving us
\begin{equation}
   \|R_1\|_1^{} \le c(s)/s^3
\end{equation}
and \eqref{q:riesz} via Young.

For $T=\agb^{-s}$, we have in place of \eqref{q:Rfkk} the stronger symmetries
\begin{equation}
   f(\pm k_1,\pm k_2)=f(k_1,k_2)
   \aand
   f(k_2,k_1)=f(k_1,k_2),
\end{equation}
resulting in better bounds for small $s$,
\begin{equation}\label{q:bdfkks}
   |\dy_i f(k_1,k_2)| \le s c(s)|k|^{-s-1} \text{ and }
   |\dy_{ij} f(k_1,k_2)| \le s c(s)|k|^{-s-2}
\end{equation}
for $i,\,j\in\{1,2\}$.
The above computations carry over almost line-by-line, except that one has to add the (subdominant) contribution of the $k_2$-axis, and that one gains a power of $s$ thanks to \eqref{q:bdfkks}.
These give \eqref{q:logb} and \eqref{q:bern1}.
\end{proof}


\nocite{rubinstein-clark:12} 
\nocite{yeung-sreenivasan:13} 

\nocite{gibson:68b} 






\bigskip\hbox to\hsize{\qquad\hrulefill\qquad}\medskip

\end{document}